\documentclass[preprint,12pt]{elsarticle}

\usepackage[left=1in, right=1in, top=1in, bottom=1in]{geometry}

\usepackage{amssymb}
\usepackage{amsmath}

\usepackage{tikz-cd}
\usepackage{graphicx}
\usepackage{amsfonts}
\usepackage{amsthm}
\usepackage{mathtools}
\usepackage{comment}
\usepackage{hyperref}

\newcommand{\nthres}[3]{\text{nthres}_{#1}^{#2}\left(#3\right)}
\newcommand{\negl}{\text{neg}}
\newcommand{\cyc}[1]{\mathbb{Z}/#1\mathbb{Z}}
\newcommand{\Hom}{\text{Hom}}

\newtheorem{theorem}{Theorem}[section]
\newtheorem{lemma}[theorem]{Lemma}
\newtheorem{corollary}[theorem]{Corollary}
\newtheorem{example}[theorem]{Example}

\theoremstyle{definition}
\newtheorem{definition}[theorem]{Definition}

\newtheorem{proposition}[theorem]{Proposition}
\newtheorem{remark}[theorem]{Remark}

\begin{document}

\begin{frontmatter}



\title{Bounds on Negligible Thresholds in Group Cohomology}


\author{Matthew Gherman, Samuel Moss} 

\affiliation{{California Institute of Technology}}

\begin{abstract}
Let $F$ be a field, $K/F$ a field extension with $\Gamma_K:=\text{Gal}(K_\text{sep}/K)$, $G$ a finite group, $j:\Gamma_K\rightarrow G$ a continuous group homomorphism, and $M$ an abelian group. There is an induced map of profinite cohomology $j^*:H^n (G;M)\rightarrow H^n (K;M)=:H^n (\Gamma_K;M)$. An element $u\in H^n (G;M)$ is negligible over $F$ if $j^* (u)=0$ for all such $j$ and field extensions $K/F$. The subset of negligible elements forms a subgroup $H^n (G;M)_{\negl,F}\subseteq H^n (G;M)$. The negligible threshold $\nthres{F}{G}{M}$ is the smallest $N\in\mathbb{N}_{\geq 0}\cup\{\infty\}$ such that $H^n (G;M)_{\negl,F}=H^n (G;M)$ for all $n\geq N$. In this paper, we first establish some general bounds for negligible thresholds of arbitrary triples $(G,M,F)$ consisting of a finite group $G$, trivial $G$-module $M$, and field $F$. We then compute the negligible thresholds of an arbitrary group when $F$ is a field of characteristic $p$ and $M=\cyc{p^k}$. Continuing in the characteristic $p$ setting, we compute these thresholds for $G$ a $p$-group and $M=\mathbb{Z}$, and we compute the thresholds for $G=\cyc{\ell},M=\cyc{s}$ for integers $s,\ell\in\mathbb{N}_{>1}$ satisfying certain number theoretic conditions. For fields containing sufficiently many roots of unity, we show that the negligible threshold of odd exponent finitely-generated trivial $G$-modules is finite. For $G$-modules of the form $\cyc{p}$, $p$ an odd prime, we compute an explicit upper bound for the negligible threshold in terms of $p$ and $|G|$. We conclude by providing explicit computations of the negligible thresholds of elementary abelian $p$-groups with mod-$p$ coefficients, negligible thresholds of cyclic groups, and bounds on mod-$p$ negligible thresholds of general finite abelian $p$-groups.
\end{abstract}





\end{frontmatter}

\tableofcontents

\section{Introduction}

The notion of negligible cohomology was introduced by J.-P. Serre in \cite{Se} (see also \cite[Part I, \S 26]{cigc}).
Let $G$ be a finite group, $M$ a $G$-module, and $F$ a field.
A continuous group homomorphism $j:\Gamma_K=\text{Gal}(K_\text{sep} /K)\to G$ from the absolute Galois group $\Gamma_K$ of a field extension $K$ of $F$ to $G$ yields a homomorphism $j^*:H^n(G;M)\to H^n(K;M)$ of cohomology groups for every $d\geq0$.
An element $u\in H^n(G;M)$ is called \textit{negligible over }$F$ if $u\in\ker(j^*)$ for all field extensions $K/F$ and all $j$.
The elements negligible over $F$ form a subgroup $H^n(G;M)_{\text{neg},F}\subseteq H^n(G;M)$.
The definition of negligible cohomology that Serre used in \cite[Part I, \S 26]{cigc} is slightly more restrictive than the one in this paper.
Serre's original definition is equivalent to our definition of negligible over $\mathbb{Q}$.
We make reference to a ground field $F$ whereas a negligible class for Serre needs to be in the kernel of the restriction map for every field.

The following examples outline important connections between negligible classes and a wide range of algebraic concepts.
\begin{example}\label{applications of negligible classes}
\begin{enumerate}
    \item[(1)] Negligible cohomology elements are related to the embedding problem.
    Let $K/F$ be a finite Galois field extension with $G=\text{Gal}(K/F)$.
    Let
    \begin{equation}\label{intro ses}
    \begin{tikzcd}
        1\arrow{r} &M\arrow{r} &G'\arrow{r}{f} &G\arrow{r} &1
    \end{tikzcd}
    \end{equation}
    be an exact sequence of finite groups with $M$ abelian.
    The conjugation $G'$-action on $M$ makes $M$ a $G$-module.
    The embedding problem for the exact sequence (\ref{intro ses}) and field extension $K/F$ is to find a Galois $G'$-algebra $K'$ over $F$ such that the restriction map $G'=\text{Gal}(K'/F)\to\text{Gal}(K/F)=G$ coincides with $f$.
    Equivalently, one needs to find a lifting $\Gamma_F\to G'$ of the homomorphism $\Gamma_F\to G$ corresponding to the extension $K/F$.
    Let $u\in H^2(G,M)$ be the class of the exact sequence (\ref{intro ses}) and let $j:\Gamma_L\to G$ be the group homomorphism given by a field extension $L/F$.
    Then $j$ extends to a homomorphism $\Gamma_L\to G'$ if and only if the pull-back of the sequence (\ref{intro ses}) under $j$ is split.
    The latter is equivalent to the triviality of the image of $u$ under $j^*:H^2(G,M)\to H^2(L,M)$.
    In other words, the class $u$ is negligible if and only if all embedding problems for the exact sequence (\ref{intro ses}) and all $G$-Galois field extensions $L'/L$ of fields containing $F$ have solutions.

    \item[(2)] Let $M$ be an abelian group which we view as a module over any profinite group with trivial action.
    The cohomology group $H^n(F;M)=H^n(\Gamma_F;M)$ is the colimit of the groups $H^n(G;M)$ over all finite discrete factor groups $G$ of $\Gamma_F$.
    The group $H^n(G;M)_{\text{neg},F}$ is contained in the kernel of the natural homomorphism $H^n(G;M)\to H^n(F,M)$.

    \item[(3)] Negligible cohomology elements of $G$ are related to the invariants of $G$ as follows.
    Let $M$ be an abelian group with trivial group action.
    Write $\text{Inv}^n(G,M)$ for the group of degree $n$ (normalized) invariants of $G$ with values in $M$ over a field $F$ (for the definition of the invariant see \cite{cigc}).
    The homomorphism
    \begin{equation*}
        \text{inv}:H^n(G;M)\to\text{Inv}^n(G,M),
    \end{equation*}
    takes $u\in H^n(G;M)$ to the invariant sending the class of a $G$-algebra $N$ over a field extension $L/F$ (that is a $G$-torsor over $\text{Spec}(L)$) to $j^*(u)\in H^n(L,M)$ for $j:\Gamma_L\to G$ the natural group homomorphism.
    By definition, $H^n(G;M)_{\text{neg},F}=\ker(\text{inv})$.
\end{enumerate}
\end{example}

The group $H^n(G;M)_{\text{neg},F}$ is trivial if $d=0$ or $d=1$ by \cite{GM}.
Under mild hypotheses for $F$, the degree two negligible cohomology was determined in \cite{MS24} for a non-trivial action on the coefficients.
The authors applied the technique to disprove a conjecture of Khare \cite{KL20}.
If $M$ is a trivial $G$-module, the degree two negligible cohomology was determined in \cite{GM}.
Serre provided explicit descriptions of mod-2 negligible classes over $\mathbb{Q}$ of elementary abelian 2-groups in \cite[Lemma 26.4]{cigc}.
The results were later extended in \cite[Theorems 4.1 and 4.3]{gherman_krull} where the mod-$p$ negligible cohomology of elementary abelian $p$-groups over arbitrary $F$ was computed.

In general, computing negligible classes is challenging.
Motivated by \cite{gherman_krull}, in which the negligible ideal in $H^*(G,\cyc{p})$ is shown to be quite large, we noticed that in many instances, when $n\gg 0$ the negligible subgroup of $H^n (G;M)$ is the entire group. We defined the negligible threshold $\nthres{F}{G}{M}$ of a triple $(G,M,F)$ (consisting of a finite group $G$, a trivial $G$-module $M$, and a field $F$) to be the minimum $N\in\mathbb{N}_{\geq 0}$ such that $H^n (G;M)_{\negl,F}=H^n (G;M)$ for every $n\geq N$. If no such minimum exists, we say $\nthres{F}{G}{M}=\infty$. 

In this paper we give bounds on these negligible thresholds for various triples $(G,M,F)$. We begin in Section \ref{2. general bounds} by establishing some elementary bounds and results on negligible thresholds. We then turn our attention to fields of characteristic $p$ in Section \ref{3. negligible threhsolds over fields of nonzero characteristic}. In particular, we provide exact threshold values in the following two cases.

\textbf{Theorem \ref{p group cyclic p}.} \textit{Let $p,k\in \mathbb{N}_{>1}$ be such that $p$ is prime. Then for any triple $(G,\cyc{p^\ell},F)$ with $\text{char}(F)=p$, $\nthres{F}{G}{\cyc{p^\ell}}=2$ if $H^1 (G;\cyc{p^\ell})\neq 0$ and $\nthres{F}{G}{\cyc{p^\ell}}=1$ otherwise.}

\textbf{Theorem \ref{p group integral}.} \textit{For any triple $(G,\mathbb{Z},F)$ with $G$ a $p$-group and $F$ a field of characteristic $p$, $\nthres{F}{G}{\mathbb{Z}}=3$.}

After this, we compute the negligible thresholds for triples of the form $(\cyc{\ell},\cyc{s},\mathbb{F}_p (\xi_s))$ when $\ell,s$ satisfy certain number theoretic constraints. In Section \ref{4. finiteness of negligible thresholds} we provide our main finiteness results.

\textbf{Theorem \ref{thresholds are often finite}.} \textit{Let $(G,M,F)$ be a triple with $M$ a finitely-generated abelian group and $F$ a field containing the roots of unity of order $|G|$ and order $\exp(M_\text{tors})$. Then $\nthres{F}{G}{M}<\infty$.}

\textbf{Theorem \ref{thresholds are often finite with a bound when p}.} \textit{Let $p$ be an odd prime and $(G,\cyc{p},F)$ a triple. Then $$\nthres{F}{G}{\cyc{p}}\leq \frac{|G|^{S+1}-1}{|G|-1}p+1.$$
Where $S:=\max\{2N,N^2\}$ and $N\in\mathbb{N}$ is the smallest positive integer such that there is a faithful $G$-representation on $\mathbb{C}^N$.}

The remainder of the paper focuses on negligible thresholds for specific groups. In Section \ref{5. thresholds mod-p of elementary abelian p-groups} we compute the thresholds for elementary abelian $p$-groups with cyclic $p$ coefficients, in Section \ref{6. negligible thresholds of cyclic groups} we compute the thresholds for cyclic groups, and in Section \ref{7. negligible thoreslds mod-p of general finite abelian p-groups} we provide threshold bounds for abelian $p$-groups with cyclic cofficients.

\subsection{Notation and definitions}

A triple $(G,M,F)$ consists of a finite group $G$, a trivial $G$-module (i.e. abelian group) $M$, and a base field $F$. We use $K /F$ to denote a field extension of $F$ and $\Gamma_K$ to denote the absolute Galois group $\text{Gal}(K_\text{sep}/K)$.

For $m\in \mathbb{N}_{\geq 2}$, we use $\zeta_m (F)$ to denote the subgroup $F^\times$ containing $m$-th roots of unity. We will write $\zeta_m$ when $F$ is clear from context. We use $\xi_m$ to denote a primitive $m$-th root of unity.

\begin{definition}
    Let $F$ be a field, $G$ a finite group, and $M$ a $G$-module. We say that an element $u\in H^n (G;M)$ is \textit{negligible over $F$} if for every field extension $K/F$ and every continuous group homomorphism $j:\Gamma_K\rightarrow G$, $j^*(u)=0 \in H^n (K;M)$. The negligible elements form a subgroup of $H^n (G;M)$, which we denote $H^n (G;M)_{\negl,F}$.
\end{definition}
\begin{definition}
    The \textit{negligible threshold} $\nthres{F}{G}{M}$ of a triple $(G,M,F)$ is the minimum $N\in \mathbb{N}_{\geq 0}\cup\{\infty\}$ for which $H^n (G;M)_{\negl,F}=H^n (G;M)$ for all $n\geq N$.
\end{definition}
From this point forward we will always assume $M$ has trivial $G$-action.

\section{General Bounds}\label{2. general bounds}

We begin with Proposition \ref{basics} providing lower bounds for triples $(G,M,F)$ with $M$ arbitrary, with $M=\mathbb{Z}$, and with $M=\mathbb{Q}/\mathbb{Z}$. In Propositions \ref{directsumsneg}, \ref{field extension} we show how the negligible threshold changes as we take direct summands of $M$ with other abelian groups, how it changes as we take extensions of $F$, and that $\nthres{F}{G}{M}\leq\nthres{\mathbb{Q}}{G}{M}$ for any triple $(G,M,F)$. Finally, in Lemma \ref{sylowthres} and Corollary \ref{sylowthres cor} we provide threshold bounds for $G$ with cylic $p$-group coefficients in terms of the Sylow $p$-subgroups of $G$.

\begin{lemma}\label{nat}
Let $G$ be a finite group, $F$ a field, and $0\rightarrow A\rightarrow B\rightarrow C\rightarrow 0$ a short exact sequence. Then the maps in the induced long exact sequence restrict to maps on the negligible subgroups.
\end{lemma}
\begin{proof}
    Let $K/F$ be a field extension with absolute Galois group $\Gamma_K$ and $\phi:\Gamma_K\rightarrow G$ a continuous group homomorphism. Then the naturality of the long exact sequence on (profinite) group cohomology yields the following diagram.
    \[\begin{tikzcd}
	{H^n (G;A)} & {H^n (G;B)} & {H^n (G;C)} & {H^{n+1}(G;C)} \\
	{H^n (K;A)} & {H^n(K;B)} & {H^n (K;C)} & {H^{n+1}(K;A)}
	\arrow[from=1-1, to=1-2]
	\arrow[from=1-1, to=2-1]
	\arrow[from=1-2, to=1-3]
	\arrow[from=1-2, to=2-2]
	\arrow[from=1-3, to=1-4]
	\arrow[from=1-3, to=2-3]
	\arrow[from=1-4, to=2-4]
	\arrow[from=2-1, to=2-2]
	\arrow[from=2-2, to=2-3]
	\arrow[from=2-3, to=2-4]
\end{tikzcd}\]
The result follows quickly from the definition of negligible classes and a simple diagram chase.

\end{proof}

\begin{lemma}\label{corestriction preserves neg}
Let $H\subseteq G$ be a subgroup of a finite group $G$.
Denote by
\begin{equation*}
    \text{cor}^d:H^n(H;M)\to H^n(G;M)
\end{equation*}
the degree $d$ corestriction map.
If $u\in H^n(H;M)$ is negligible over $F$, then $\text{cor}^d(u)\in H^n(G;M)_{\text{neg},F}$.
\end{lemma}

\begin{proof}
Let $K=F(V)^G$ for $V$ a faithful $F$-representation of $G$.
Since $H$ is a subgroup of $G$, $V$ is likewise a faithful $F$-representation of $H$.
Define $K_H=F(V)^H$.
We build the following commutative square.
\begin{equation*}
\begin{tikzcd}
    H^n(H;M)\arrow{r}{\text{cor}^d}\arrow{d}{\text{inf}_H} &H^n(G;M)\arrow{d}{\text{inf}} \\
    H^n(K_H;M)\arrow{r}{\text{cor}^d} &H^n(K;M)
\end{tikzcd}
\end{equation*}
By \cite[Proposition 2.1]{GM}, the kernel of these inflation maps are the negligible classes of $H^n(H;M)$ and $H^n(G;M)$, respectively.
Then $\text{inf}_H(u)=0$ implies that $\text{inf}(\text{cor}^d(u))=0$ as desired.

\end{proof}

\begin{proposition}\label{basics}
Let $(G,M,F)$ be a triple. 

(1) If $H^1 (G;M)\neq 0$, then $\nthres{F}{G}{M}\geq 2$.

(2) If $G$ is not perfect, then $\nthres{F}{G}{\mathbb{Z}}\geq 3$. 

(3) If $\nthres{F}{G}{\mathbb{Z}}\geq 3$, then $\nthres{F}{G}{\mathbb{Q}/\mathbb{Z}}=\nthres{F}{G}{\mathbb{Z}}-1$. 
\end{proposition}

\begin{proof}
  By \cite[Corollary 2.2]{GM}, $H^n (G;M)_\negl =0$ for $i=0,1$, thus $H^1 (G;M)\neq 0$ implies $H^1 (G;M)_\negl\neq H^1 (G;M)$. This proves (1).

Since $\mathbb{Q}$ is an injective $\mathbb{Z}[G]$-module, $H^n (G;\mathbb{Q})=0$ for $n>0$, so the short exact sequence $0\rightarrow\mathbb{Z}\rightarrow\mathbb{Q}\rightarrow\mathbb{Q}/\mathbb{Z}\rightarrow 0$ implies that the connecting homomorphism $\partial:H^n (G;\mathbb{Q}/\mathbb{Z})\rightarrow H^{n+1}(G;\mathbb{Z})$ is an isomorphism. By Lemma \ref{nat}, it follows that $H^n (G;\mathbb{Q}/\mathbb{Z})_{\negl,F}\cong H^{n+1} (G;\mathbb{Z})_{\negl,F}$ when $n>1$. 

There are well-known isomorphisms $$H^1 (G;\mathbb{Q}/\mathbb{Z})\cong \Hom_{\mathsf{Grp}}(G,\mathbb{Q}/\mathbb{Z})\cong \Hom_\mathsf{Ab} (G/[G,G],\mathbb{Q}/\mathbb{Z}).$$
Suppose $G$ is not perfect. Then $G\neq [G,G]$, so $G/[G,G]$ is a finitely-generated abelian group. By the classification theorem, there exist $r,\ell_1,\dots,\ell_r,p_1,\dots,p_r\in\mathbb{N}$ with the $p_i$ prime such that $G/[G,G]\cong \bigoplus_{i=1}^r \cyc{p_i^{\ell_k}}$. Therefore,
$$\Hom_\mathsf{Ab}(G/[G,G],\mathbb{Q}/\mathbb{Z})\cong \Hom_\mathsf{Ab}(\bigoplus_{i=1}^r \cyc{p_i^{\ell_i}},\mathbb{Q}/\mathbb{Z})\cong\bigoplus_{i=1}^r \Hom_\mathsf{Ab}(\cyc{p_i^{\ell_i}},\mathbb{Q}/\mathbb{Z}).$$
But for each $i$ there is an embedding $\cyc{p_i^{\ell_i}}\hookrightarrow\mathbb{Q}/\mathbb{Z}$, so the homset above is nontrivial. By \cite[Corollary 2.2(d)]{GM}, we see that $0=H^1 (G;\mathbb{Q}/\mathbb{Z})_{\negl,F}\subsetneq H^1 (G;\mathbb{Q}/\mathbb{Z})$, so by the first paragraph $0=H^2 (G;\mathbb{Z})_{\negl,F}\subsetneq H^2 (G;\mathbb{Z})$, which proves (2).

Claim (3) follows from the isomorphisms $H^n (G;\mathbb{Q}/\mathbb{Z})_{\negl,F}\cong H^{n+1}(G;\mathbb{Z})_{\negl,F}$ and $H^n (G;\mathbb{Q}/\mathbb{Z})\cong H^{n+1}(G;\mathbb{Z})$ when $n>1$.  

\end{proof}

The first part of Example \ref{cyclic char p example} shows that the bound in (1) is sharp.
Further, Theorem \ref{p group integral} shows the bound in (2) is sharp.

\begin{proposition}\label{relprimeneg}
    Let $(G,\mathbb{Z}/\ell\mathbb{Z},F)$ be a triple with $\ell\in\mathbb{N}_{>1}$ and $\gcd(|G|,\ell)=1$. Then $\nthres{F}{G}{\mathbb{Z}/\ell\mathbb{Z}}=1$.
\end{proposition}
\begin{proof}
    Since $\gcd(|G|,\ell)=1$, $H^n (G;\cyc{\ell})=0$ for $n>0$.
    Since $\cyc{\ell}$ has a trivial $G$-action, $H^0 (G;\cyc{\ell})=\cyc{\ell}$. Thus $\nthres{F}{G}{\cyc{\ell}}=1$.
    
\end{proof}

\begin{proposition}\label{directsumsneg}

     Let $\{(G,M_i,F):i\in I\}$ be a collection of triples. Then $$\nthres{F}{G}{\bigoplus_{i\in I}M_i}=\sup_{i\in I}\nthres{F}{G}{M_i}.$$
\end{proposition}
\begin{proof}
    For any extension $K/F$ and morphism $\phi:\Gamma_K\rightarrow G$, the following commutes for all $n\in\mathbb{N}$

\begin{center}
\begin{tikzcd}
	{H^n (G;\bigoplus_{i\in I}M_i)} & {\bigoplus_{i\in I}H^n (G;M_i)} \\
	{H^n (K;\bigoplus_{i\in I} M_i)} & {\bigoplus_{i\in I}H^n (K;M_i)}
	\arrow["\cong", from=1-1, to=1-2]
	\arrow["{\phi^*_{\bigoplus_{i\in I}M_i}}"', from=1-1, to=2-1]
	\arrow["{\bigoplus_{i\in I}\phi^*_{M_i}}", from=1-2, to=2-2]
	\arrow["\cong"', from=2-1, to=2-2]
\end{tikzcd}
\end{center}
If $n<\sup_{i\in I}\nthres{F}{G}{M_i}$, then there exists some $j\in I$, extension $K'/F$, morphism $\phi':\Gamma_{K'}\rightarrow G$, and $u_j\in H^n (G;M_j)$ such that $\phi^*_{M_j}(u_j)\neq 0$, therefore $\bigoplus_{i\in I}\phi^*_{M_i}(u_i)_{i\in I}=(\phi^*_{M_i} (u_i))_{i\in I}\neq 0$, so by the diagram above $\phi^*_{\bigoplus_{i\in I}M_i} (u)\neq 0$, where $u$ is image of $(u_i)_{i\in I}$ under the bottom isomorphism, therefore $\nthres{F}{G}{\bigoplus_{i\in I}M_i}\geq \sup_{i\in I}\nthres{F}{G}{M_i}$.

On the other hand, if $n>\sup_{i\in I}\nthres{F}{G}{M_i}$, then for all extensions $K/F$ and morphisms $\phi:\Gamma_K\rightarrow G$, $\phi^*_{M_i}=0$ for all $i\in I$, thus the r.h.s. map in the diagram is zero, thus so is the left one, so $\nthres{F}{G}{\bigoplus_{i\in I}M_i}\leq\sup_{i\in I}\nthres{F}{G}{M_i}$.

\end{proof}

\begin{proposition}\label{field extension}
    Let $(G,M,F)$ be a triple and $K/F$ a field extension. Then 
    
    (1) $\nthres{K}{G}{M}\leq\nthres{F}{G}{M}$.

    (2) For any triple $(G,M,F)$, $\nthres{F}{G}{M}\leq\nthres{\mathbb{Q}}{G}{M}$.
\end{proposition}

\begin{proof}
    Let $n\geq\nthres{F}{G}{M}$. By \cite[Proposition 2.3(2)]{GM}, $H^n (G;M)_{\negl,F}\subseteq H^n (G;M)_{\negl,K}$, therefore
$$H^n(G;M)_{\negl,F}=H^n (G;M)\subseteq H^n (G;M)_{\negl,K}\subseteq H^n (G;M)$$
Therefore $H^n (G;M)_{\negl,K}=H^n (G;M)$, which proves (1).

For (2), let $n\geq\nthres{\mathbb{Q}}{G}{M}$. The proof of \cite[Proposition 26.1]{cigc} shows $H^n (G;M)_{\negl,\mathbb{Q}}\subseteq H^n (G;M)_{\negl,F}$.
Therefore, $$H^n (G;M)=H^n (G;M)_{\negl,\mathbb{Q}}\subseteq H^n (G;M)_{\negl,F}\subseteq H^n (G;M).$$
Thus, $H^n (G;M)_{\negl,F}=H^n (G;M)$.

\end{proof}

The bound in (1) is clearly sharp by taking $K=F$, but Corollary \ref{bound for Z/sZ} together with Example \ref{cyclic char p example} reveal we can have equality for a non-trivial field extension.
The bound in (2) is clearly sharp by taking $F=\mathbb{Q}$, but \cite[Theorem 4.3]{gherman_krull} proves that the mod $p>2$ negligible classes of an elementary abelian $p$-group are the same over $\mathbb{Q}$ and any field extension $F$ that does not contain a primitive $p^2$-th root of unity.

\begin{proposition}\label{ordergneg}
    For a finite group $G$ and field $F$, $\nthres{F}{G}{\mathbb{Z}}\leq\nthres{F}{G}{\cyc{|G|}}+1$
\end{proposition}

\begin{proof}
    From properties of restriction-corestriction, the map $m_{|G|}:H^n (G;\mathbb{Z})\rightarrow H^n (G;\mathbb{Z})$ is the zero map for all $n\in\mathbb{N}_{\geq 0}$. From the short exact sequence $0\rightarrow\mathbb{Z}\rightarrow\mathbb{Z}\rightarrow\cyc{|G|}\rightarrow 0$ we have for each $n\in\mathbb{N}_{\geq 0}$ exact sequences
\[\begin{tikzcd}
	{H^n (G;\mathbb{Z}/|G|\mathbb{Z})_{\text{neg},F}} & {H^{n+1}(G;\mathbb{Z})_{\text{neg},F}} & {H^{n+1}(G;\mathbb{Z})_{\text{neg},F}}
	\arrow["{\partial}",from=1-1, to=1-2]
	\arrow["{m_{|G|}}", from=1-2, to=1-3]
\end{tikzcd}\]
By exactness $\text{im }\partial=\ker m_{|G|}=H^{n+1}(G;\mathbb{Z})$. From the naturality of long exact sequences and the definition of negligible cohomology, this restricts to a map $\partial:H^n (G;\cyc{|G|})_{\negl,F}\rightarrow H^{n+1} (G;\mathbb{Z})_{\negl,F}$. Now take $n\geq \nthres{F}{G}{\cyc{|G|}}$. Then $H^n (G;\cyc{|G|})_{\negl,F}=H^n (G;\cyc{|G|})$, so $\text{im }\partial=H^{n+1}(G;\mathbb{Z})\subseteq H^{n+1} (G;\mathbb{Z})_{\negl,F}$, which proves the claim.

\end{proof}

Theorems \ref{p group cyclic p} and \ref{p group integral} show the bound in Proposition \ref{ordergneg} is sharp for $G$ a $p$-group with $F$ characteristic $p$.

\begin{lemma}\label{sylowthres}
Let $(G,\cyc{p^\ell},F)$ be a triple with $p$ prime. Then for any Sylow-$p$ subgroup $H$ of $G$,
$$\nthres{F}{G}{\cyc{p^\ell}}\leq \nthres{F}{H}{\cyc{p^\ell}}$$
\end{lemma}

\begin{proof}
We know that $\text{cor}\circ\text{res}:H^n (G;\cyc{p^\ell})\rightarrow H^n (H;\cyc{p^\ell})\rightarrow H^n (G;\cyc{p^\ell})$ is multiplication by $[G:H]$. Since $H^n (G;\cyc{p^\ell})$ is finite and $p$-torsion, this must be an isomorphism, thus $\text{cor}$ is surjective. By Lemma \ref{corestriction preserves neg}, the following diagram commutes 
\[\begin{tikzcd}
	{H^n (H;\cyc{p^\ell})_{\negl,F}} & {H^n (H;\cyc{p^\ell})} \\
	{H^n (G;\cyc{p^\ell})_{\negl,F}} & {H^n (G;\cyc{p^\ell})}
	\arrow[hook, from=1-1, to=1-2]
	\arrow["{\text{cor}}"', two heads, from=1-1, to=2-1]
	\arrow["{\text{cor}}", two heads, from=1-2, to=2-2]
	\arrow[hook, from=2-1, to=2-2]
\end{tikzcd}\]
If $H^n (H;\cyc{p^\ell})_{\negl,F}=H^n (H;\cyc{p^\ell})$, a simple diagram chase shows $H^n (G;\cyc{p^\ell})_{\negl,F}=H^n (G;\cyc{p^\ell})$.

\end{proof}

\begin{corollary}\label{sylowthres cor}
    Let $(G,M,F)$ be a triple with $M=\bigoplus_{i=1}^s \cyc{p^{\ell_i}_i}$ a finite abelian group. Then
    $$ \nthres{F}{G}{M}\leq \max_{i\in \{1,\dots,s\}, H_{p_i}\in \text{Syl}_{p_i}(G)}\nthres{F}{H_{p_i}}{\cyc{p^{\ell_i}_i}} $$
\end{corollary}
\begin{proof}
    This follows immediately from Proposition \ref{directsumsneg} and Lemma \ref{sylowthres}.
    
\end{proof}

Note that Lemma \ref{sylowthres} and Corollary \ref{sylowthres cor} are sharp when $G$ is a $p$-group and $M=\mathbb{Z}/p^\ell\mathbb{Z}$.

\section{Negligible Thresholds over Fields of Nonzero Characteristic}\label{3. negligible threhsolds over fields of nonzero characteristic}

When $F$ is a field of characteristic $p$, the negligible thresholds of many triples $(G,M,F)$ are at most 3. We first show that in Theorem \ref{p group cyclic p} that if $M$ is a cyclic $p$-group, this negligible threshold is either 1 or 2. We then show that if $M=\mathbb{Z}$ and $G$ is a $p$-group, the negligible threshold is 3 in Theorem \ref{p group integral}. Finally, we restrict our attention to $G$ and $M$ cyclic groups with orders satisfying certain number-theoretic conditions in Propositions \ref{trivial K group} and \ref{bound for Z/sZ}, where we show that the negligible threshold is either 2 or 3, and that both thresholds can occur.

\begin{theorem}\label{p group cyclic p}
Let $p,k\in \mathbb{N}_{>1}$ be such that $p$ is prime. Then for any triple $(G,\cyc{p^\ell},F)$ with $\text{char}(F)=p$, $\nthres{F}{G}{\cyc{p^\ell}}=2$ if $H^1 (G;\cyc{p^\ell})\neq 0$ and $\nthres{F}{G}{\cyc{p^\ell}}=1$ otherwise.
\end{theorem}

\begin{proof}
    By the Artin-Schreier sequence and Hilbert's Theorem 90 \cite[6.11.16(1)]{weibelhomalg}, $H^n (L;\cyc{p})=0$ for all $n> 1$. Using the short exact sequences
$$0\rightarrow \cyc{p}\rightarrow\cyc{p^2}\rightarrow\cyc{p}\rightarrow0$$
$$0\rightarrow\cyc{p}^{k-1}\rightarrow\cyc{p}^k\rightarrow\cyc{p}\rightarrow0$$
and proceeding by induction, it can be easily shown that $H^n (L;\cyc{p}^k)=0$ for all $n\geq 1$.

Let $K/F$ be any field extension and $\phi:\Gamma_K\rightarrow G$ any morphism. Then by the argument above, $H^n (K;\cyc{p}^k)=0$ for any $n,k\in\mathbb{N}_{>1}$, so $\phi^*_{\cyc{p^k}}:H^n (G;\cyc{p^k})\rightarrow H^n (K;\cyc{p^k})=0$ is the zero map for any such $n,k$. Thus $\nthres{F}{G}{\cyc{p^k}}\leq 2$.  The result follows from Proposition \ref{basics}(1) and \cite[Corollary 2.2(2)]{GM} that $H^1 (G;\cyc{p}^k)_\negl=0$.

\end{proof}

\begin{theorem}\label{p group integral}
    For any triple $(G,\mathbb{Z},F)$ with $G$ a $p$-group and $F$ a field of characteristic $p$, $\nthres{F}{G}{\mathbb{Z}}=3$.
\end{theorem}

\begin{proof}
    Let $|G|=p^\ell$. Consider the short exact sequence $0\rightarrow\mathbb{Z}\rightarrow\mathbb{Z}\rightarrow\cyc{p^\ell}\rightarrow 0$, where the first map is multiplication by $p^\ell$. This induces a map $m_{p^\ell}:H^n (G;\mathbb{Z})\rightarrow H^n (G;\mathbb{Z})$ which sends $u\mapsto p^\ell\cdot u$.
    Cohomology on the short exact sequence above yields an exact sequence $$H^{n-1}(G;\cyc{p^\ell})\rightarrow H^n (G;\mathbb{Z})\xrightarrow{m_{p^\ell}}H^n (G;\mathbb{Z}).$$
    It is well-know that $m_{p^\ell}=0$, so the leftmost arrow is surjective.
    By Theorem \ref{p group cyclic p}, $H^{n-1} (G;\cyc{p^\ell})_\negl=H^{n-1} (G;\cyc{p^\ell})$ for $n\geq3$, and since this map restricts to one on negligible cohomology, we see that $H^n (G;\mathbb{Z})_\negl=H^n (G;\mathbb{Z})$ for $n\geq 3$.
    Thus, $\nthres{F}{G}{\mathbb{Z}}\leq 3$.
    Since a $p$-group is never perfect, Proposition \ref{basics}(2) shows $\nthres{F}{G}{\mathbb{Z}}\geq 3$.
    
\end{proof}

In particular, taking Proposition \ref{directsumsneg}, Theorem \ref{p group cyclic p}, and Theorem \ref{p group integral} together, we get a complete computation of negligible thresholds of triples $(G,M,F)$ with $G$ a $p$-group, $F$ of characteristic $p$, and $M$ a finitely-generated abelian group.

By completely different methods, we can compute the negligible thresholds for certain cyclic groups over fields of characteristic $p$.

\begin{proposition}\label{trivial K group}
Let $p\in\mathbb{N}$ be prime and $s,\ell\in\mathbb{N}$ such that $\gcd(p,s)=1$ and $\ell | s$. Then $\nthres{\mathbb{F}_p (\xi_s)}{\cyc{\ell}}{\cyc{s}}\in \{2,3\}$.
\end{proposition}

\begin{proof}
Let $L:=\mathbb{F}_p (\xi_s)(t^\ell)$. By the Norm Residue Isomorphism Theorem \cite{HW}, for any $n\in\mathbb{N}$ we have $H^n (L;\cyc{s})\cong K_n^M (L)/s K_n^M (L)$, where $K_n^M (L)$ is the $n$-th Milnor $K$-group of $L$. Since $L$ has transcendental degree 1 over $\mathbb{F}_p (\xi_s)$, it is a global field, so by \cite[Section 2]{BassTate1973Milnor} $K_n^M (L)=0$ whenever $n\geq 3$. Thus $H^n (L;\cyc{s})=0$ for $n\geq 3$.

Since $\ell$ divides $s$ and $\gcd(p,s)=1$, there is a faithful action of $\cyc{\ell}$ on $\mathbb{F}_p (\xi_s)(t)$ generated by $t\mapsto \zeta_s^{s/\ell} t$. It is clear that the fixed subfield under this action is exactly $L$. By the argument above, the canonical map $j:\Gamma_L\rightarrow \cyc{\ell}$ induces the zero map on cohomology of degree at least 3. The fact that $\nthres{\mathbb{F}_p}{\cyc{\ell}}{\cyc{s}}\leq 3$ follows from \cite[Proposition 2.1]{GM}.

Since $\ell$ divides $s$, there is a well-defined nonzero morphism $\cyc{\ell}\rightarrow\cyc{s}$ sending $[1]\mapsto [s/\ell]$, therefore $H^1 (\cyc{\ell};\cyc{s})=\Hom_\mathsf{Grp} (\cyc{\ell},\cyc{s})\neq 0$. By Proposition \ref{basics}(1), $\nthres{\mathbb{F}_p}{\cyc{\ell}}{\cyc{s}}\geq 2$.

\end{proof}

Keeping the same notation as above, let $|p|_s^\times$ denote the order of the class $[p]\in(\cyc{s})^\times$.

\begin{proposition}\label{bound for Z/sZ}
If $\gcd(|p|_s^\times,s)=1$, then $\nthres{\mathbb{F}_p}{\cyc{\ell}}{\cyc{s}}\in\{2,3\}$. More generally, if $F$ has characteristic $p$ then $\nthres{F}{\cyc{\ell}}{\cyc{s}}\in\{2,3\}$.
\end{proposition}

\begin{proof}
From the theory of finite fields, $[\mathbb{F}_p (\xi_s):\mathbb{F}_p]=\min\{k\in\mathbb{Z}:s \mid p^k-1\}$. But note that $s \mid p^k-1$ iff there exists $n\in\mathbb{Z}$ such that $s n=p^k-1$ iff $p^k\equiv 1 (\text{mod } s)$, thus the smallest such $k$ is precisely $|p|_s^\times$. By \cite[Proposition 2.3(3)]{GM}, $[\mathbb{F}_p (\xi_s):\mathbb{F}_p]\cdot H^n (\cyc{\ell};\cyc{s})_{\negl,\mathbb{F}_p (\xi_s)}\subseteq H^n (\cyc{\ell};\cyc{s})_{\negl,\mathbb{F}_p}$. The computation of cohomology of cyclic groups implies $|H^n (\cyc{\ell};\cyc{s})|$ divides $s$. But $[\mathbb{F}_p (\xi_s):\mathbb{F}_p]=|p|_s^\times$ is relatively prime to $s$ by assumption, thus multiplication of $H^n (\cyc{\ell};\cyc{s})$ by $[\mathbb{F}_p (\xi_s):\mathbb{F}_p]$ is an isomorphism. Since this map clearly respects negligible classes, it follows that $H^n (\cyc{\ell};\cyc{s})_{\negl,\mathbb{F}_p (\xi_s)}=H^n (\cyc{\ell};\cyc{s})_{\negl,\mathbb{F}_p}$. In particular, when $n\geq \nthres{\mathbb{F}_p (\xi_s)}{\cyc{\ell}}{\cyc{s}}$ we have $H^n (\cyc{\ell};\cyc{s})_{\negl,\mathbb{F}_p (\xi_s)}=H^n (\cyc{\ell};\cyc{s})$, so $H^n (\cyc{\ell};\cyc{s})_{\negl,\mathbb{F}_p}=H^n (\cyc{\ell};\cyc{s})$. Thus $\nthres{\mathbb{F}_p}{\cyc{\ell}}{\cyc{s}}\leq \nthres{\mathbb{F}_p (\xi_s)}{\cyc{\ell}}{\cyc{s}}$. From Proposition \ref{field extension}(1) and the fact that characteristic $p$ fields are extensions of $\mathbb{F}_p$, $\nthres{F}{\cyc{\ell}}{\cyc{s}}\leq \nthres{\mathbb{F}_p}{\cyc{\ell}}{\cyc{s}}\leq\nthres{\mathbb{F}_p (\xi_s)}{\cyc{\ell}}{\cyc{s}}$. By Proposition \ref{basics}(1) and the fact that $H^1 (\cyc{\ell};\cyc{s})\neq 0$, $2\leq \nthres{F}{\cyc{\ell}}{\cyc{s}}$, proving both claims.

\end{proof}

Given a field $F$, let $\zeta_s$ denote the $s$-th roots of unity in some algebraic closure $\overline{F}$ of $F$.
With the setup of Proposition \ref{trivial K group}, \cite[Proposition 3.1]{GM} and the isomorphism $\zeta_s\cong\cyc{s}$ imply$$H^2 (\cyc{\ell};\cyc{s})_{\negl,\mathbb{F}_p (\xi_s)}=\text{im}(H^1 (\cyc{\ell};\mathbb{F}_p (\xi_s)^\times/\zeta_s)\rightarrow H^2 (\cyc{\ell};\zeta_s))$$
where the r.h.s. map comes from the short exact sequence $1\rightarrow\zeta_s\rightarrow\mathbb{F}_p (\xi_s)^\times\rightarrow \mathbb{F}_p (\xi_s)^\times/\zeta_s\rightarrow 1$.
Thus, to determine whether the threshold is 2 or 3, we just need to compute $H^2 (\cyc{\ell};\cyc{s})_{\negl,\mathbb{F}_p (\xi_s)}$.
If it is equal to the group cohomology $H^2 (\cyc{\ell};\cyc{s})$, then the threshold is 2, otherwise it is 3.
Example \ref{cyclic char p example} illustrates the method.

\begin{example}\label{cyclic char p example}
Consider the case $s=\ell=3,p=19$. Since $3\cdot 6=18$ and $\mathbb{F}_{19}^\times$ is cyclic of order 18, it follows that $\zeta_3\subseteq \mathbb{F}_p^\times$, so $\mathbb{F}_{19}(\zeta_3)=\mathbb{F}_{19}$.
Since the $\cyc{3}$-action is trivial, the induced map $H^1 (\cyc{3};\zeta_3)\rightarrow H^1 (\cyc{3};\mathbb{F}_{19} (\zeta_3)^\times)$ is injective.
With these observations and the first isomorphism theorem in mind, the formula from earlier implies there is an exact sequence
$$0\rightarrow H^1 (\cyc{3};\cyc{3})\rightarrow H^1 (\cyc{3};\cyc{18})\rightarrow H^1 (\cyc{3};\cyc{6})\rightarrow H^2 (\cyc{3};\cyc{3})_{\negl,\mathbb{F}_{19}}\rightarrow 0.$$
Since the first cohomology is the group homset and since the homsets of $\cyc{3}$ with $\cyc{3},\cyc{18},\cyc{6}$ are all isomorphic to $\cyc{3}$, the sequence above becomes
$$0\rightarrow\cyc{3}\rightarrow\cyc{3}\rightarrow\cyc{3}\rightarrow H^2 (\cyc{3};\cyc{3})_{\negl,\mathbb{F}_{19}}\rightarrow 0.$$
By exactness, $H^2 (\cyc{3};\cyc{3})_{\negl,\mathbb{F}_{19}}\neq 0$.
From the computation of cyclic group cohomology,
\begin{equation*}
    H^2 (\cyc{3};\cyc{3})\cong\cyc{3}
\end{equation*}
and, thus, is a simple group, so it must be the case that $H^2 (\cyc{3};\cyc{3})_{\negl,\mathbb{F}_{19}}=H^2 (\cyc{3};\cyc{3})$.
This proves $\nthres{\mathbb{F}_{19}}{\cyc{3}}{\cyc{3}}=2$.

Now take $s=9,\ell\in \{3,9\},p=19$. By similar reasoning as before, $\mathbb{F}_{19}(\zeta_9)=\mathbb{F}_{19}$. Note that $|\mathbb{F}_{19}^\times/\zeta_9|=2$ is relative prime to $3$ and $9$, thus $H^1 (\cyc{\ell};\mathbb{F}_{19}^\times / \zeta_9)=0$.
By the formula above, $H^2 (\cyc{\ell};\cyc{9})_{\negl,\mathbb{F}_{19}}=0$.
But cyclic group cohomology is 2-periodic, thus $H^2 (\cyc{\ell};\cyc{9})\cong H^0 (\cyc{\ell};\cyc{9})=\cyc{9}$.
This implies  $H^2 (\cyc{\ell};\cyc{9})_{\negl,\mathbb{F}_{19}}\subsetneq H^2 (\cyc{\ell};\cyc{9})$, so $\nthres{\mathbb{F}_{19}}{\cyc{9}}{\cyc{\ell}}=3$.
\end{example}

\section{Finiteness of Negligible Thresholds}\label{4. finiteness of negligible thresholds}

In this section, we first show in Lemma \ref{cycfin} that the negligible threshold is finite for triples $(G,M,F)$ such that $M$ is an odd cyclic group and $F$ contains all $|M|$.
In Theorem \ref{thresholds are often finite}, we generalize the result to any triple $(G,M,F)$ in which $M$ is a finitely-generated abelian group and $F$ contains sufficient roots of unity.
In the case that $M$ has odd prime order, we provide an explicit upper bound in Theorem \ref{thresholds are often finite with a bound when p} for the negligible threshold and argue that it is likely not sharp. Finally, we mention that, if $G$ has even order and $M=\mathbb{Z}/2\mathbb{Z}$, then the negligible threshold could be infinite.

\begin{lemma}\label{cycfin}
    Fix odd $\ell\in\mathbb{N}_{>1}$ and let $(G,\cyc{\ell},F)$ be a triple with $\zeta_\ell \subseteq F$. Then $\nthres{F}{G}{\cyc{\ell}}<\infty$.
\end{lemma}

\begin{proof}

Let $F$ be a field, $K/F$ any extension. Let $n$ be odd and take $\alpha\in H^n (K;\cyc{\ell})$. Since $H^* (K;\cyc{\ell})$ is anticommutative, $\alpha\cup \alpha=(-1)^{\deg(\alpha)\deg(\alpha)}\alpha\cup\alpha=-\alpha\cup\alpha$. Since $\ell$ is odd, this implies $\alpha\cup\alpha=0$. Let $H^\text{even}(K;\cyc{\ell})\subseteq H^* (K;\cyc{\ell})$ be the subalgebra containing only the even degree pieces. By anticommutativity of $H^* (K;\cyc{\ell})$, this is a commutative ring. By the Norm Residue Isomorphism Theorem \cite{HW} and the fact that $\zeta_\ell\subseteq F\subseteq K$, $H^*(K;\cyc{\ell})$ is generated as an $H^0 (K;\cyc{\ell})$-algebra by $H^1 (K;\cyc{\ell})$. By the argument above, any element of the form $u_1\cup\dots\cup u_m$ for $u_i\in H^1 (K;\cyc{\ell})$ squares to zero. As the $H^0 (K;\cyc{\ell})$-span of the irrelevant ideal of $H^\text{even}(K;\cyc{\ell})$ and since the nilradical of a commutative ring is an ideal, we have that all elements of $H^n (K;\cyc{\ell})$ are nilpotent for $n>0$.

By \cite[Proposition 2.1]{GM}, it suffices to show that for large enough $n$, $H^n (G;\cyc{\ell})=\ker j^*$ for $j:\Gamma_{F(V)^G}\rightarrow G$ the canonical map. It is shown in \cite[Corollary 1.14]{hoegenhaven2013} that $H^* (G;\cyc{\ell})$ is finitely-generated as a $\cyc{\ell}$-algebra. Decomposing each of these generators into their homogeneous pieces, we may assume without loss of generality that the generators are homogeneous. Let $\{u_1,\dots,u_m\}$ be such a generating set. By the previous paragraph, for $1\leq i\leq m$ there exist $r_i\in\mathbb{N}$ such that $j^* (u_i)^{r_i}=0$. Let $N:=m\max_{1\leq i\leq m}r_i$. Take $n>N$ and $u\in H^n (G;\cyc{\ell})$. Since $\{u_1,\dots,u_m\}$ is a generating set, $u=\sum_{i=1}^r k_i \prod_{j=1}^m u_j^{s_{i,j}}$ for some $r,s_{i,j}\in\mathbb{N},k_i\in\cyc{\ell}$ with $\sum_{j=1}^m s_{i,j}=n$. Since $j^*$ is a map of algebras,
$ j^* (u)=\sum_{i=1}^r k_i\prod_{j=1}^m j^*(u_j)^{s_{i,j}}$. By the pigeonhole principle, for some $j$ we have $s_{i,j}>r_j$, thus $j^* (u_j)^{s_{i,j}}=0$, so $j^* (u)=0$. Thus $u\in H^n (G;\cyc{\ell})_{\negl,F}$.

\end{proof}

\begin{corollary}\label{intfin}
Let $(G,\mathbb{Z},F)$ be a triple with $\zeta_{|G|}\subseteq F$. Then $\nthres{F}{G}{\mathbb{Z}}<\infty$.
\end{corollary}
\begin{proof}
    By Proposition \ref{ordergneg}, $\nthres{F}{G}{\mathbb{Z}}$ differs from $\nthres{F}{G}{\cyc{|G|}}$ by at most 1. Thus the result follows from Lemma \ref{cycfin}.
    
\end{proof}

\begin{theorem}\label{thresholds are often finite}
Let $(G,M,F)$ be a triple with $M$ a finitely-generated abelian group and $F$ a field containing the roots of unity of order $|G|$ and order $\exp(M_\text{tors})$. Then $\nthres{F}{G}{M}<\infty$.
\end{theorem}

\begin{proof}
Since $M$ is finitely-generated, $M\cong\mathbb{Z}^r\oplus \bigoplus_{i=1}^s \cyc{p^{n_i}_i}$ for $p_1,\dots,p_s$ primes and $n_1,\dots,n_s\in \mathbb{N}_{>0}$. By Proposition \ref{directsumsneg},
$$ \nthres{F}{G}{M}=\max_{i\in \{1,\dots,s\}}\{\nthres{F}{G}{\mathbb{Z}},\nthres{F}{G}{\cyc{p^{n_i}_i}}\} $$
Since we have assumed sufficiently many roots of unity, Lemma \ref{cycfin} and Corollary \ref{intfin} imply $\nthres{F}{G}{\mathbb{Z}}<\infty$ and $\nthres{F}{G}{\cyc{p^{n_i}_i}}<\infty$ for each $i$, hence the result.

\end{proof}

\begin{theorem}\label{thresholds are often finite with a bound when p}
   Let $p$ be an odd prime and $(G,\cyc{p},F)$ a triple. Then $$\nthres{F}{G}{\cyc{p}}\leq \frac{|G|^{S+1}-1}{|G|-1}p+1$$
Where $S:=\max\{2N,N^2\}$ and $N\in\mathbb{N}$ is the smallest positive integer such that there is a faithful $G$-representation on $\mathbb{C}^N$.  
\end{theorem}

\begin{proof}

By \cite[Corollary 2.4]{GM}, it suffices to show this bound holds for $F(\zeta_p)$. Let $N$ be the smallest integer such that there is a faithful $G$-representation on $\mathbb{C}^N$. By \cite[Corollary 4.3]{totaro_group_2014}, the maximum degree of the generators of $H^* (G;\cyc{p})$ is $S:=\max\{2N,N^2\}$. Using the Bar resolution $B_\bullet$, the cohomology of the following chain complex yields the group cohomology
$$0\rightarrow\Hom_{\mathbb{Z}[G]}(\mathbb{Z},\cyc{p})\rightarrow\Hom_{\mathbb{Z}[G]}(B_1,\cyc{p})\rightarrow\Hom_{\mathbb{Z}[G]}(B_2,\cyc{p})\rightarrow\Hom_{\mathbb{Z}[G]}(B_3,\cyc{p})\rightarrow\dots$$
From the construction of the Bar resolution, $B_n$ is free with $|G|^n$ generators. Now any $\mathbb{Z}[G]$-morphism $\phi:B_n\rightarrow\cyc{p}$ is determined by where the generators of $B_n$ gets mapped to. For each of the $|G|^n$ generators there are $p$ possibilities, thus $|\Hom_{\mathbb{Z}[G]}(B_n,\cyc{p})|=p^{|G|^n}$, so $\dim_{\cyc{p}}\Hom_{\mathbb{Z}[G]} (B_n,\cyc{p})=|G|^n$. Since the cohomology $H^n (G;\cyc{p})$ is a quotient of a subvector space of $\Hom_{\mathbb{Z}[G]}(B_n,\cyc{p})$, we have $\dim_{\cyc{p}} H^n (G;\cyc{p})\leq |G|^n$. Thus there are at most $\sum_{i=1}^S |G|^n=\frac{|G|^{S+1}-1}{|G|-1}$ generators for $H^* (G;\cyc{p})$.

Now let $j:\Gamma_{F(V)^G}\rightarrow G$ be the canonical map. From the proof of Lemma \ref{cycfin}, $u\cup u=0$ for any $u \in H^1 (F(V)^G;\cyc{p})$. By the Norm Residue Isomorphism Theorem \cite{HW}, any $v\in H^* (F(V)^G;\cyc{p})$ is of the form $v=\sum_{i=1}^r k_i u_i$ for $r\in\mathbb{N},k_i\in\cyc{p}$, and $u_i\in H^1 (F(V)^G;\cyc{p})$. By freshman's dream and since $u_i^p=u_i^2 u_i^{p-2}=0$, $v^p=\sum_{i=1}^r k_i^p u_i^p=0$. Thus $j^* (w)^p=0$ for any $w\in H^*(G;\cyc{p})$. The claim then follows from using these bounds in the proof of Lemma \ref{cycfin}.

\end{proof}

\begin{remark}
    This bound is in general not sharp. Consider the triple $(G,M,F):=(\cyc{3},\cyc{3},\mathbb{F}_3)$. As shown in Theorem \ref{p group cyclic p}, $\nthres{\mathbb{F}_3}{\cyc{3}}{\cyc{3}}=2$. We know that there is a faithful representation of $\cyc{3}$ on $\mathbb{C}$ sending the generator to a primitive root of unity, so in the notation of Theorem 4.2, $N=1$ and $S=2$. Then Theorem \ref{thresholds are often finite with a bound when p} gives the bound $\frac{3^3-1}{3-1}\cdot3+1=40$. It is interesting to note, however, that this bound does not depend on the field $F$.
\end{remark} 

\begin{remark}
    There are cases where the threshold is infinite. Let $G$ have even order and $F$ be formally real. Note that the set of all $F$-negligible elements of the commutative ring $H^*(G;\cyc{2})$ form an ideal. By \cite[Theorem 3.3]{gherman_krull}, the quotient by this ideal has Krull dimension 1. If the negligible threshold were finite, this quotient must be a finite ring and thus have Krull dimension 0. We conclude that $\nthres{F}{G}{\cyc{2}}=\infty$.
\end{remark}








\section{Negligible Thresholds Mod-$p$ of Elementary Abelian $p$-Groups}\label{5. thresholds mod-p of elementary abelian p-groups}

In this section, we compute bounds on the negligible thresholds of elementary abelian $p$-groups over various fields with coefficients $\cyc{p}$.
We make use of the computations of the negligible ideal given in \cite[Section 4]{gherman_krull}. We begin with the case $p=2$.

\subsection{The case $p=2$}

Let $G=(\cyc{2})^n$ for some $n\in\mathbb{N}_{>0}$. Then
\begin{equation*}
    H^* (G;\cyc{2})=\cyc{2}[x_1,\dots,x_n]
\end{equation*}
with $\deg x_i=1$ for each $1\leq i\leq n$. If $F$ is formally real, then $\nthres{F}{G}{\cyc{2}}=\infty$ by \cite[Theorem 3.3]{gherman_krull}. Let $F$ be a field with characteristic not equal to 2 and which is not formally real. Let $s(F)=2^r$ denote the level of $F$.
By \cite[Section 2]{elman2008algebraic}, $(-1)^{r+1}\in H^*(F,\mathbb{Z}/2\mathbb{Z})$ is trivial.


\begin{theorem}\label{elementary abelian 2-group thresh}
Let $G$ be an elementary abelian 2-group of rank $n$. Let $F$ be a field which is not formally real with $s(F)=2^r$ for some $r\in\mathbb{N}_{\geq 0}$. Then $\nthres{F}{G}{\cyc{2}}=n(r+1)+1$.
\end{theorem}

\begin{proof}
    We split into two cases: $r=0$ and $r>0$. First consider the case $r>0$, so that $s(F)=2^r>1$. By \cite[Theorem 4.1(b)]{gherman_krull}, $I(G,\cyc{2})_F$ is generated by $\{x_i x_j^2+x_j x_i^2:1\leq i<j\leq n\}\cup \{x_i^{r+2}:1\leq i\leq n\}$. Let $m\geq n(r+1)+1$. Then any $H^m (G;\cyc{2})$ is generated by $\{x_1^{r_1}\dots x_n^{r_n}:r_i\in\mathbb{N}_{\geq 0},\sum_{i=1}^n r_i=m \}$. By the pigeonhole principle, for some $1\leq j\leq n$ we have $r_j>r+1$, thus each term making up $u$ contains an $x_j^{r+2}$ part and so is negligible. Thus $u$ is negligible, and $\nthres{F}{G}{\cyc{2}}\leq n(r+1)+1$. 

Take $x_1^{r+1}\dots x_n^{r+1}\in H^{n(r+1)}(G;\cyc{2})$. Since the degree of each $x_i$ is less than $r+2$, in order to show this is negligible it suffices to show that it is not contained in the ideal generated by $\{x_i x_j^2+x_j x_i^2:1\leq i<j\leq n\}$. Suppose this is not the case. Then there exist $p_{i,j}\in \cyc{2}[x_1,\dots,x_n]$ such that $$x_1^{r+1}\dots x_n^{r+1}=\sum_{1\leq i<j\leq n} p_{i,j}(x_1,\dots,x_n)(x_i x_j^2+x_j x_i^2)$$
Evaluating at $x_i=1$ for each $i$, the $x_i x_j^2+x_j x_i^2$ map to $1+1=0$, thus the r.h.s. vanishes. But $x_1^{r+1}\dots x_n^{r+1}$ evaluates to $1$, so the equation above implies $0\neq 1$, which is absurd. Thus $x_1^{r+1}\dots x_n^{r+1}\notin I(G;\cyc{2})_F$, so $\nthres{F}{G}{\cyc{2}}>n(r+1)$, which proves the claim when $s(F)=2^r>1$.

Now suppose $r=0$, so that $s(F)=1$. By \cite[Theorem 4.1(c)]{gherman_krull}, $I(G;\cyc{2})_F$ is generated by $\{x_i^2:1\leq i\leq n\}$. Since $\{x_1^{r_1}\dots x_n^{r_n}:r_i\in\mathbb{N}_{\geq 0},\sum_{i=1}^n r_i=m\}$ generates $H^m(G;\cyc{2})$, it follows from a similar argument as the first paragraph that when $m\geq n+1$, each generator is negligible, thus $H^m(G;\cyc{2})$ is entirely negligible, so $\nthres{F}{G}{\cyc{2}}\leq n+1$. On the other hand, it is clear that $x_1\dots x_n$ is not in the ideal generated by the $\{x_i^2:1\leq i\leq n\}$, so $H^n (G;\cyc{2})$ is not entirely negligible, and $\nthres{F}{G}{\cyc{2}}>n$.

\end{proof}

\subsection{The odd prime case}

For $p\in\mathbb{N}$ an odd prime and $n\in\mathbb{N}_{>1}$, let $G:=(\cyc{p})^n$. By \cite[Section 3.5]{Ev}, the cohomology ring is given by $H^* (G;\cyc{p})=\Lambda_{\cyc{p}}(G^*)[y_1,\dots,y_n]$, where $\deg y_i=2$ and $\Lambda_{\cyc{p}}(G^*)$ is the exterior algebra of the character group $H^1 (G;\cyc{p})$ of $G$. Let $\{x_1,\dots,x_n\}$ be a basis for $H^1 (G;\cyc{p})$. By the description of the cohomology ring given above, $x_i\cup x_j=0$ iff $i=j$.

\begin{theorem}\label{elementary abelian threshold}
If $F$ contains a primitive $p^2$-th root of unity, then $\nthres{F}{G}{\cyc{p}}= n+1$. If $F$ does not contain a primitive $p^2$-th root of unity, then $\nthres{F}{G}{\cyc{p}}= n+2$.
\end{theorem}

\begin{proof}
Suppose $F$ contains a primitive $p^2$-th root of unity. By \cite[Theorem 5.3(b)]{gherman_krull}, the negligible ideal $I(G;\cyc{p})_F$ is generated by $\{y_i:1\leq i\leq n\}$. The degree one part is generated by some $x_1,\dots,x_n$. Clearly $x_1\cup \dots\cup x_n$ does not lie in this ideal, thus $H^n (G;\cyc{p})$ has a non-negligible element, so $\nthres{F}{G}{\cyc{p}}>n$. Theorem \ref{finite abelian odd p}(3) implies the other inequality.

Now suppose $F$ does not contain a primitive $p^2$-th root of unity. By \cite[Theorem 4.3(a)]{gherman_krull}, the negligible ideal is generated by $\{x_i y_j+x_j y_i:1\leq i\leq j\leq n\}\cup\{y_i y_j:1\leq i\leq j\leq n\}$. Take $u\in H^m (G;\cyc{p})$ for $m\geq n+3$. Then each term of $u$ must contain at least two $y_i$ terms. Since the ring is anticommutative, it follows that each term of $u$ (and thus $u$) has a $y_i y_j$ for some $i,j$, thus $u$ is negligible. It follows that $\nthres{F}{G}{\cyc{p}}\leq n+3$. 

Now consider a monomial of degree $n+2$. Since each $x_i$ can only appear once in a monomial, this monomial either has at least two $y_i$ or is of the form $x_1\cup\dots\cup x_n\cup y_i$ for some $1\leq i\leq n$. In the former case, it is immediate the element is negligible. In the latter case, note that up to a $\pm 1$ is equal to $ x_1\cup\dots\cup\hat{x}_i\cup\dots\cup x_n\cup x_i\cup y_i=\frac{1}{2}x_1\cup\dots\cup\hat{x}_i\cup\dots\cup x_n \cup(x_i y_i+x_i y_i)$, where $\hat{x}_i$ means this term is excluded and $1/2$ denotes the multiplicative inverse of two in $\cyc{p}$ (which exists since $p$ is odd). Thus elements of this form is also negligible. As these two cases form a $\cyc{p}$-basis for $H^{n+2}(G;\cyc{p})$, it follows that $H^{n+2}(G;\cyc{p})_{\negl,F}=H^{n+2}(G;\cyc{p})$, thus $\nthres{F}{G}{\cyc{p}}\leq n+2$.

Consider the element $x_1\cup\dots\cup x_{n-1}\cup y_n\in H^{n+1} (G;\cyc{p})$. Suppose this lies in the negligible ideal. Then for $1\leq i\leq j\leq n$ there exist $p_{ij},q_{ij}\in H^* (G;\cyc{p})$ such that
$$x_1\cup\dots\cup x_{n-1}\cup y_n=\sum_{1\leq i\leq j\leq n} (p_{ij}\cdot (x_i y_j+x_j y_i)+q_{ij} \cdot y_i y_j)$$

We may assume without loss of generality that the $p_{ij}$ contain no terms with a $y_i$, since this could be absorbed into the corresponding $q_{i'j'}$.
Since the left hand side has $y$-degree 1, the $q_{ij}$ must be zero.
Thus
$$x_1\cup\dots\cup x_{n-1}\cup y_n=\sum_{1\leq i\leq j\leq n} p_{ij}\cdot (x_i y_j+x_j y_i),\hspace{5mm}p_{ij}\in\Lambda_{\cyc{p}} (G^*)$$
Note that by definition $p_{ij}=p_{ji}$. We will show by induction that
$$ \sum_{1\leq i\leq j\leq n}p_{ij}\cdot (x_i y_j+x_j y_i)=\sum_{i=1}^n (2p_{ii}x_i+\sum_{j=1,j\neq i}^n p_{ij}x_j)y_i$$
If $n=1$, this is immediate. Now suppose the result holds for $n$. By the induction hypothesis and $p_{ij}=p_{ji}$,
\begin{align*}
    \sum_{1\leq i\leq j\leq n+1} p_{ij}\cdot (x_i y_j+x_j y_i)&=\sum_{1\leq i\leq j\leq n} p_{ij}\cdot (x_i y_j+x_j y_i)+\sum_{i=1}^{n+1} p_{i,n+1}(x_i y_{n+1}+x_{n+1}y_i)\\ &=\sum_{i=1}^n (2 p_{ii}x_i+\sum_{j=1,j\neq i}^n p_{ij}x_j)y_i+\sum_{i=1}^{n}p_{i,n+1}x_i y_{n+1}\\ &\hspace{10mm}+\sum_{i=1}^n p_{i,n+1}x_{n+1}y_i+2 p_{n+1,n+1}x_{n+1}y_{n+1}\\ &=\sum_{i=1}^n (2p_{ii}x_i+\sum_{j=1,j\neq i}^n p_{ij}x_j+p_{i,n+1}x_{n+1})y_i+(2 p_{n+1,n+1}+\sum_{i=1}^n p_{i,n+1}x_i)y_{n+1}\\ &=\sum_{i=1}^{n+1}(2p_{ii}x_i+\sum_{j=1,j\neq i}^{n+1}p_{ij}x_j)y_i
\end{align*}
This completes the induction step. Since the left hand side and the $x_iy_j+x_j y_i$ are homogeneous of degress $n+1$ and $3$ (resp.), we may assume that each term in the $p_{ij}$ is homogeneous of degree $n-2$. Since $\Lambda_{\cyc{p}}(G^*)$ is generated by $x_1,\dots,x_n$ (each of degree 1) and since $x_i^2=0$, we must have
$p_{ij}=\sum_{s,t\in [n], s\neq t}\alpha^{ij}_{st}x_{[n]\setminus \{s,t\}}$, where $\alpha^{ij}_{st}\in\cyc{p}$, $[n]=\{1,\dots,n\}$, and $x_{\{i_1<i_2<\dots<i_\ell\}}=x_{i_1}x_{i_2}\dots x_{i_\ell}$. Note that $\alpha_{st}^{ij}=\alpha^{ij}_{ts}$. Then
$p_{ij}x_k=\sum_{s,t\in [n],s\neq t}\alpha^{ij}_{st}x_{[n]\setminus \{s,t\}}x_k$. When $s,t\neq k$, $x_{[n]\setminus \{s,t\}}x_k$ has an $x_k^2$ term, thus is zero. It follows that the sum above is equal to $\sum_{s\in [n]\setminus \{k\}}\alpha^{ij}_{sk}x_{[n]\setminus\{s,k\}}x_k$. For $s\neq k$, define
$$\epsilon^n_{sk}:=\begin{cases}
    (-1)^{n-k} &\text{if } k>s\\ (-1)^{n-k+1} &\text{if }k<s
\end{cases}$$
Note that $x_{[n]\setminus \{s,k\}}x_k=\epsilon^n_{sk}x_{[n]\setminus \{s\}}$. With all this in mind, the $y_i$ coefficient of $\sum_{1\leq i\leq j\leq n}p_{ij}(x_i y_j+x_j y_i)$ can be rewritten:
$$\sum_{i=1}^n \left(2\sum_{s\in [n]\setminus \{i\}}\alpha^{ii}_{si}\epsilon^n_{si}x_{[n]\setminus\{s\}}+\sum_{j=1,j\neq i}^n \sum_{s\in [n]\setminus \{j\}}\alpha^{ij}_{sj}\epsilon^n_{sj}x_{[n]\setminus \{s\}}\right)$$
By linear independence of the $x_{[n]\setminus \{s\}}$ across $s\in [n]$, the equation above implies the $x_{[n]\setminus \{i\}}$ coefficient of the $y_i$ coefficient is given by
$$\sum_{j=1,j\neq i}^n \alpha^{ij}_{ij}\epsilon_{ij}^n=\sum_{j=1}^{i-1}(-1)^{n-j+1}\alpha^{ij}_{ij}+\sum_{j=i+1}^n (-1)^{n-j}\alpha^{ij}_{ij}$$
For ease of notation we will write $\beta_{ij}:=\alpha_{ij}^{ij}$; note that $\beta_{ij}=\beta_{ji}$. We will show by induction that
$$\sum_{i=1}^n (-1)^i \left(\sum_{j=1}^{i-1}(-1)^{n-j+1}\beta_{ij}+\sum_{j=i+1}^n (-1)^{n-j}\beta_{ij}\right)=0\hspace{10mm}(*)$$
This is tautologically true when $n=1$. Now suppose it holds for $n$. Then for $n+1$ we have by the induction hypothesis and $\beta_{ij}=\beta_{ji}$ that
\begin{align*}
    \sum_{i=1}^{n+1} (-1)^i \big(\sum_{j=1}^{i-1}(-1)^{(n+1)-j+1}\beta_{ij}&+\sum_{j=i+1}^{n+1} (-1)^{(n+1)-j}\beta_{ij}\big)\\ &=(-1)\sum_{i=1}^n (-1)^i \left(\sum_{j=1}^{i-1}(-1)^{n-j+1}\beta_{ij}+\sum_{j=i+1}^{n}(-1)^{n-j}\beta_{ij}\right)\\ &\hspace{10mm}+\sum_{i=1}^n (-1)^i (-1)^{n-(n+1)+1}\beta_{i,n+1}+(-1)^{n+1}\sum_{j=1}^n (-1)^{n-j+2}\beta_{n+1,j}\\ &=\sum_{i=1}^n (-1)^i\beta_{i,n+1}+\sum_{j=1}^n (-1)^{2n-j+3}\beta_{j,n+1}\\ &=\sum_{i=1}^n (-1)^i \beta_{i,n+1}-\sum_{i=1}^n (-1)^i\beta_{i,n+1}=0
\end{align*}
Where the last line follows from the fact that $2n-i+3\equiv i+1\hspace{.5mm}(\text{mod }2)$. This completes the induction. But note that since $\sum_{1\leq i\leq j\leq n}p_{ij}\cdot (x_i y_j+x_j y_i)=x_1\cup\dots\cup x_{n-1}\cup y_n$, the $x_{[n]\setminus \{i\}}$ coefficient of the $y_i$ coefficient is $0$ for all $i\neq n$ and is $1$ for $i=n$. Thus $(*)$ must evaluate to $(-1)^n$: a contradiction! It follows that the original assumption is wrong, thus $x_1\cup\dots\cup x_{n-1}\cup y_n$ is not contained in the negligible ideal, so $\nthres{F}{G}{\cyc{p}}>n+1$. 

\end{proof}

\section{Negligible Thresholds of Cyclic Groups}\label{6. negligible thresholds of cyclic groups}

By the end of the section, we will prove Theorem \ref{nthresh cyclic group} about negligible thresholds by computing the negligible classes in degrees three and four for cyclic groups with some cyclic coefficients.

\begin{theorem}\label{nthresh cyclic group}
Let $G$ be cyclic of order $m$ and $F$ a field with $\text{char}(F)$ relatively prime to $m$.
If $m\not\equiv 2$ (mod 4), assume $\zeta_m\subseteq F$.
If $m\equiv 2$ (mod 4), assume $\zeta_{2m}\subseteq F$.
Then $\nthres{F}{G}{\mathbb{Z}/\ell\mathbb{Z}}\in\{2,3\}$ for $\ell|m$ and $\nthres{F}{G}{\mathbb{Z}}=3$.
\end{theorem}

\begin{remark}\label{sharp cyclic threshold}
To determine whether $\nthres{F}{G}{\mathbb{Z}/m\mathbb{Z}}$ is 2 or 3 for cyclic $G$, we only need to apply \cite[Corollary 3.4]{GM}.
The relevant factor is the roots of unity present in the base field $F$.
In particular, $\nthres{F}{G}{\mathbb{Z}/m\mathbb{Z}}=3$ when $G$ is a cyclic group of odd order $m$, $\zeta_{m}\subseteq F$, and $\zeta_{m^2}\not\subseteq F$ and $\nthres{F}{G}{\mathbb{Z}/m\mathbb{Z}}=2$ when $G$ is a cyclic group of odd order $m$ and $\zeta_{m^2}\subseteq F$
\end{remark}

We will need preliminary results about connecting homomorphisms in order to prove Theorem \ref{nthresh cyclic group}.
Let $K$ be a field extension of $F$ and $j:\Gamma_K\to G$ a continuous group homomorphism.
We will assume throughout this section that $\text{char}(F)$ is relatively prime to $m$ and $\zeta_m\subseteq F$, so $\zeta_m$ is a trivial $\Gamma_K$-module.
Endow $\text{Ext}^1_{\Gamma_K}(\zeta_m,\zeta_m)$ with a group structure via the Baer sum.
The group homomorphism 
\begin{equation*}
    \Phi_k:\text{Ext}^1_{\Gamma_K}(\zeta_m,\zeta_m)\to\text{Hom}(H^k(K;\zeta_m),H^{k+1}(K;\zeta_m))
\end{equation*}
identifies the class of an extension with a connecting homomorphism in the induced long exact sequence on cohomology.

Let $C_i$ in $\text{Ext}^1_{\Gamma_K}(\zeta_m,\zeta_m)$ denote the class of the extension
\begin{equation*}
\begin{tikzcd}
    1\arrow{r} &\zeta_m\arrow{r} &\zeta_{m^2}^{\otimes i}\arrow{r} &\zeta_m\arrow{r} &1.
\end{tikzcd}
\end{equation*}
Let $\xi_{m^2}$ be a primitive $m^2$ root of unity.
The cyclic extension $K(\xi_{m^2})/K$ provides a group homomorphism $\phi:\Gamma_K\to\text{Gal}(K(\xi_{m^2})/K)$.
The action on $\zeta_{m^2}^{\otimes i}$ is given by $\sigma\cdot \xi_m^\ell=\phi(\sigma)^i(\xi_m^\ell)$.

For a cyclic degree $m$ field extension $E$ of $K$, we know that $E=K(\sqrt[m]{x})$ for some $x\in E^\times$ since $\zeta_m\subseteq K$.
Denote by $\chi:\Gamma_K\to\text{Gal}(E/K)\cong\mathbb{Z}/m\mathbb{Z}$ the surjective group homomorphism induced by the field extension.
Let $D(x)$ in $\text{Ext}^1_{\Gamma_K}(\zeta_m,\zeta_m)$ be the class of the extension
\begin{equation*}
\begin{tikzcd}
    1\arrow{r} &\zeta_m\arrow{r} &\zeta_m\oplus\zeta_m\arrow{r} &\zeta_m\arrow{r} &1.
\end{tikzcd}
\end{equation*}
with action $\sigma\cdot(\xi_1,\xi_2)=(\xi_1\chi(\sigma)(\xi_2),\xi_2)$ for $\sigma\in\Gamma_K$.

\begin{lemma}\label{Baer sum}
Let $K$ be a field and $\xi_m$ a primitive $m$-th root of unity.
Assume $\zeta_m\subseteq K$.
Then the Baer sum $kD(\xi_m)+C_0$ is equal to $C_k$ in $\text{Ext}^1_{\Gamma_K}(\zeta_m,\zeta_m)$.
\end{lemma}

\begin{proof}
Suppose $\zeta_{m^2}\subseteq K$.
Then the sequence for $D(\xi_{m^2})$ splits, so this is class is zero.
Furthermore, $\zeta^{\otimes k}_{m^2}=\zeta_{m^2}$, so $C_k=C_0$.
Thus, $k D(\xi_m)+C_0=C_k$ holds trivially.

Now assume $\zeta_{m^2}\not\subseteq K$. We will proceed by induction on $k$.
The base case $k=0$ is trivial.
Assume that $kD(\xi_m)+C_0=C_k$.
We will show that the Baer sum $D(\xi_m)+C_k=C_{k+1}$.
Let $X$ be the pullback in the following diagram.
\begin{equation*}
\begin{tikzcd}
    \zeta_m\arrow[hookrightarrow]{r}\arrow[dashed]{rd} &\zeta_m\oplus\zeta_m\arrow[twoheadrightarrow]{r} &\zeta_m \\
    &X\arrow{u}\arrow{r} &\zeta_{m^2}^{\otimes k}\arrow[twoheadrightarrow]{u} \\
    & &\zeta_m\arrow[hookrightarrow]{u}\arrow[dashed]{ul}
\end{tikzcd}
\end{equation*}
Then $X=\{((\xi,\xi'),\xi''):\xi,\xi'\in\zeta_m,\xi''\in\zeta_{m^2},\xi'=(\xi'')^m\}$.
Define $Y$ to be the quotient of $X$ via the image of the skew-diagonal embedding of $\zeta_m$ in $X$.
The image of the skew-diagonal embedding of $\zeta_m$ in $X$ is $\{((\xi^{-1},1),\xi):\xi\in\zeta_m\}$ so each class of $Y$ has a representative of the form $((1,\xi'),\xi'')$ for $\xi'\in\zeta_m$, $\xi''\in\zeta_{m^2}$, and $\xi'=(\xi'')^m$. Moreover, we see the representatives of this form are uniquely determined by $\xi''$, thus the abelian group homomorphism $Y\rightarrow\zeta_{m^2}^{\otimes (k+1)}$ sending $((1,\xi'),\xi'')\mapsto\xi''$ is a bijection.

The Baer sum $D(\xi_m)+C_k$ is the class of the extension $0\to\zeta_m\to Y\to\zeta_m\to 0$.

Both $D(\xi_m)$ and $C_k$ refer to the cyclic field extension $K(\xi_{m^2})$ of $K$.
Denote the group homomorphism corresponding to the extension by $\chi:\Gamma_K\to\text{Gal}(K(\xi_{m^2})/K)$.
By \cite[Theorem VI.6.2(b)]{lang_algebra_2005} and basic Galois theory, there exists $\sigma\in\Gamma_K$ such that $\chi(\sigma)$ generates $\text{Gal}(K(\xi_{m^2})/K)$ and $\chi(\sigma)$ is multiplication by $(\xi'')^m$ for $\xi''\in\zeta_{m^2}$ a primitive $m^2$-root of unity.
In other words, the action of $\sigma$ is multiplication by $(\xi'')^m$.
Then
\begin{align*}
    \sigma\cdot((1,\xi'),\xi'')&=((\chi(\sigma)(\xi'),\xi'),\chi(\sigma)^k(\xi'')) \\
    &=((\xi',\xi'),\chi(\sigma)^k(\xi'')) \\
    &=((1,\xi'),\chi(\sigma)^k(\xi'')(\xi')) \\
    &=((1,\xi'),\chi(\sigma)^{k+1}(\xi''))
\end{align*}
since $\xi'=(\xi'')^m$. Since $((1,\xi'),\xi'')$ generates $Y$, we conclude $Y\cong\zeta_{m^2}^{\otimes(k+1)}$ as $\Gamma_K$-modules and $D(\xi_m)+C_k=C_{k+1}$ in $\text{Ext}^1_{\Gamma_K}(\zeta_m,\zeta_m)$.

\end{proof}

\begin{lemma}\label{ext computation}
Let $L$ be a field with $\zeta_m\subseteq L$.
Then $\Phi_k(C_k)$ is trivial and $\Phi_k(C_0)=-k\Phi_k(D(\xi_m))$ in $\Hom(H^k(L;\zeta_m),H^{k+1}(L;\zeta_m))$.
\end{lemma}

\begin{proof}
It is sufficient to show $\Phi_k(C_k)$ is trivial in $\Hom(H^k(L;\zeta_m),H^{k+1}(L;\zeta_m))$ by Lemma \ref{Baer sum}.
Since $\zeta_m\subseteq L$, we can identify $\zeta_m\cong\zeta_m^{\otimes i}$ as $\Gamma_L$-modules for all $i\in\mathbb{Z}$.
The extension 
\begin{equation*}
\begin{tikzcd}
    1\arrow{r} &\zeta_m\arrow{r} &\zeta_{m^2}^{\otimes d}\arrow{r} &\zeta_m\arrow{r} &1
\end{tikzcd}
\end{equation*}
induces a long exact sequence in cohomology with portion
\begin{equation*}
\begin{tikzcd}
    H^n(L;\zeta_{m^2}^{\otimes d})\arrow{r} &H^n(L;\zeta_m)\arrow{r}{\partial^n} &H^{n+1}(L;\zeta_m).
\end{tikzcd}
\end{equation*}
By the Norm Residue Isomorphism Theorem \cite{HW}, the following commutative square
\begin{equation*}
\begin{tikzcd}
    H^n(L;\zeta_{m^2}^{\otimes n})\arrow{r} &H^n(L;\zeta_m^{\otimes n}) \\
    K_n^M(L)/m^2K_n^M(L)\arrow{u}{\cong}\arrow[twoheadrightarrow]{r} &K_n^M(L)/mK_n^M(L)\arrow{u}{\cong}
\end{tikzcd}
\end{equation*}
proves that the top map is surjective.
Therefore, $\Phi_k(C_k)=\partial^k$ is trivial.

\end{proof}

\begin{lemma}
Let $K$ be a field that contains a primitive $m$-th root of unity $\xi_m$.
Then
\begin{equation*}
    \Phi_k(D(\xi_m))(\gamma)=(\xi_m)\cup\gamma
\end{equation*}
for $\gamma\in H^k(K;\zeta_m)$.
\end{lemma}

\begin{proof}
Let $\xi_{m^2}$ be a primitive $m^2$-th root of unity.
For the cyclic field extension $K(\xi_{m^2})$ of $K$, there is a corresponding group homomorphism $\chi:\Gamma_K\to\text{Gal}(K(\xi_{m^2})/K)\hookrightarrow \mathbb{Z}/m\mathbb{Z}$ in $H^1(K;\mathbb{Z}/m\mathbb{Z})$ corresponding to the generator of the Galois group.
The pairing $\mathbb{Z}/m\mathbb{Z}\otimes\zeta_m\to\zeta_m$ induces a cup product
\begin{equation*}
    \cup:H^1(K;\mathbb{Z}/m\mathbb{Z})\otimes H^k(K;\zeta_m)\to H^{k+1}(K;\zeta_m)
\end{equation*}
for each $k\geq0$.
Since $\zeta_m\subseteq K$, we can identify $D(\xi_m)$ with the short exact sequence 
\begin{equation*}
\begin{tikzcd}
    0\arrow{r} &\mathbb{Z}/m\mathbb{Z}\arrow{r} &\mathbb{Z}/m\mathbb{Z}\oplus\mathbb{Z}/m\mathbb{Z}\arrow{r} &\mathbb{Z}/m\mathbb{Z}\arrow{r} &0
\end{tikzcd}
\end{equation*}
with $\Gamma_K$-action on $\mathbb{Z}/m\mathbb{Z}\oplus\mathbb{Z}/m\mathbb{Z}$ given by $\sigma\cdot(x,y)=(x+\chi(\sigma)(y),y)$ for $\sigma\in\Gamma_K$.
Then $D(\xi_m)$ induces the commutative square
\begin{equation*}
\begin{tikzcd}
    H^0(K;\mathbb{Z}/m\mathbb{Z})\otimes H^k(K;\zeta_m)\arrow{r}{\cup}\arrow{d}{\partial^0\otimes\text{id}} &H^k(K;\zeta_m)\arrow{d}{\partial^k} \\
    H^1(K;\mathbb{Z}/m\mathbb{Z})\otimes H^k(K;\zeta_m)\arrow{r}{\cup} &H^{k+1}(K;\zeta_m)
\end{tikzcd}
\end{equation*}
for $\partial^k:H^k(K;\zeta_m)\to H^{k+1}(K;\zeta_m)$ a connecting map.
We identify $H^0(K;\mathbb{Z}/m\mathbb{Z})\cong\mathbb{Z}/m\mathbb{Z}$ so $\partial^0$ maps $1\in\mathbb{Z}/m\mathbb{Z}$ to $\chi\in H^1(K;\mathbb{Z}/m\mathbb{Z})$.
By commutativity,
\begin{equation*}
    \Phi_k(D(\xi_m))(\gamma)=\partial^k(\gamma)=\chi\cup\gamma
\end{equation*}
for $\gamma\in H^k(K;\zeta_m)$.
In the notation we adopt, $\chi$ is written as $(\xi_m)\in K^\times/(K^\times)^m$.

\end{proof}

\begin{lemma}\label{Z to Z/m surjection}
Let $K$ be a field that contains a primitive $m$-th root of unity $\xi_m$.
The mod-$m$ surjection $q:\mathbb{Z}\to\mathbb{Z}/m\mathbb{Z}$ induces a map $q_*:H^2(K;\mathbb{Z})\to H^2(K;\mathbb{Z}/m\mathbb{Z})$ on cohomology.
A character $\chi\in H^2(K;\mathbb{Z})$ corresponds to a cyclic extension $K(\sqrt[m]{a})$ of $K$ defined by $q\circ\chi:\Gamma_K\to\mathbb{Z}/m\mathbb{Z}$.
Then $q_*(\chi)=(a)\cup(\xi_m)$.
\end{lemma}

\begin{proof}
The following commutative diagram
\begin{equation*}
\begin{tikzcd}
    0\arrow{r} &\mathbb{Z}\arrow{r} &\mathbb{Q}\arrow{r} &\mathbb{Q}/\mathbb{Z}\arrow{r} &0 \\
    0\arrow{r} &\mathbb{Z}\arrow{r}\arrow[equal]{u}\arrow{d}{q} &\mathbb{Z}\arrow{r}\arrow{u}{\frac{1}{m}}\arrow{d} &\mathbb{Z}/m\mathbb{Z}\arrow{r}\arrow[hookrightarrow]{u}\arrow[equal]{d} &0 \\
    0\arrow{r} &\mathbb{Z}/m\mathbb{Z}\arrow{r} &\mathbb{Z}/m^2\mathbb{Z}\arrow{r} &\mathbb{Z}/m\mathbb{Z}\arrow{r} &0
\end{tikzcd}
\end{equation*}
induces the commutative square
\begin{equation*}
\begin{tikzcd}
    H^1(K;\mathbb{Q}/\mathbb{Z})\arrow{r}{\cong} &H^2(K;\mathbb{Z})\arrow{d}{q_*} \\
    H^1(K;\mathbb{Z}/m\mathbb{Z})\arrow{u}\arrow{r}{\partial^1} &H^2(K;\mathbb{Z}/m\mathbb{Z}).
\end{tikzcd}
\end{equation*}
The image of $(a)$ in $H^2(K;\mathbb{Z})$ is $\chi$.
By Lemma \ref{ext computation},
\begin{equation*}
    \partial^1(a)=-\Phi_1(C_0)=-(\xi_m)\cup(a)=(a)\cup(\xi_m).
\end{equation*}
The result follows from commutativity.

\end{proof}

\begin{proposition}\label{cyclic Z/mZ coefficients}
Let $G$ be a cyclic group of order $m$.
If $m\not\equiv 2$ (mod 4), assume $\zeta_m\subseteq F$.
If $m\equiv 2$ (mod 4), assume $\zeta_{2m}\subseteq F$.
Then $H^n(G;\mathbb{Z}/m\mathbb{Z})_{\text{neg},F}=H^n(G;\mathbb{Z}/m\mathbb{Z})$ for $d\in\{3,4\}$.
\end{proposition}

\begin{proof}
Let $K$ be a field extension of $F$ and $j:\Gamma_K\to G$ a continuous group homomorphism.
Let $x\in H^1(G;\mathbb{Z}/m\mathbb{Z})$ be a generator.
Then $j^*(x)\in H^1(K;\mathbb{Z}/m\mathbb{Z})$ corresponds to a field extensions $E/K$ for which $\text{Gal}(E/K)\cong G\cong\mathbb{Z}/m\mathbb{Z}$.
Since $\zeta_m\subseteq F$, we note $E=K(\sqrt[m]{a})$ for some $(a)\in K^\times/(K^\times)^m$.
Then $j^*(x)$ corresponds to $(a)$ in $H^1(K;\mathbb{Z}/m\mathbb{Z})\cong K^\times/(K^\times)^m$.
By including $\mathbb{Z}/m\mathbb{Z}$ into $\mathbb{Q}/\mathbb{Z}$, the field extension induces a character $\chi_G:G\to\mathbb{Q}/\mathbb{Z}$ in $H^1(G;\mathbb{Q}/\mathbb{Z})\cong H^2(G;\mathbb{Z})$.

Let $q:\mathbb{Z}\to\mathbb{Z}/m\mathbb{Z}$ be a surjection with induced map $q_*:H^2(-;\mathbb{Z})\to H^2(-;\mathbb{Z}/m\mathbb{Z})$ on cohomology.
Denote $\chi=j^*(\chi_G)\in H^2(K;\mathbb{Z})$ and generator $y=q_*(\chi_G)\in H^2(G;\mathbb{Z}/m\mathbb{Z})$.
By Lemma \ref{Z to Z/m surjection}, the following commutative square proves $j^*(y)=(a)\cup(\xi_m)\in H^2(K;\mathbb{Z}/m\mathbb{Z})$.
\begin{equation*}
\begin{tikzcd}
    H^2(G;\mathbb{Z})\arrow{r}{q_*}\arrow{d}{j^*} &H^2(G;\mathbb{Z}/m\mathbb{Z})\arrow{d}{j^*} \\
    H^2(K;\mathbb{Z})\arrow{r}{q_*} &H^2(K;\mathbb{Z}/m\mathbb{Z})
\end{tikzcd}
\end{equation*}
Since the cohomology of cyclic groups is 2-periodic, the cup product by $\chi_G\in H^2(G;\mathbb{Z})$ gives isomorphisms $H^n(G;\mathbb{Z}/m\mathbb{Z})\to H^{n+2}(G;\mathbb{Z}/m\mathbb{Z})$ in each degree $n\geq 1$.
Thus $x\cup y$ is a generator of $H^3(G;\mathbb{Z}/m\mathbb{Z})$ and $y\cup y$ is a generator of $H^4(G;\mathbb{Z}/m\mathbb{Z})$.
We have
\begin{align*}
    j^*(x\cup y)&=(a)\cup(a,\xi_m)=(a,-1,\xi_m) \\
    j^*(y\cup y)&=(a,\xi_m)\cup(a,\xi_m)=-(a,-1,\xi_m,-1).
\end{align*}
When $m\not\equiv 2$ (mod 4) and $\zeta_m\subseteq F$, we note $(-1)$ is trivial.
When $m\equiv 2$ (mod 4), we have
\begin{align*}
    j^*(x\cup y)&=(a,-1,\xi_m)=(a,-1,\xi_{2m}^2)=2(a,-1,\xi_{2m})=0 \\
    j^*(y\cup y)&=-(a,-1,\xi_m,-1)=-(a,-1,\xi_{2m}^2,-1)=-2(a,-1,\xi_{2m},-1)=0.
\end{align*}
Therefore, $H^3(G;\mathbb{Z}/m\mathbb{Z})$ and $H^4(G;\mathbb{Z}/m\mathbb{Z})$ are entirely negligible.

\end{proof}

\begin{corollary}\label{cyclic Z coefficients}
Let $G$ be cyclic of order $m$.
If $m\not\equiv 2$ (mod 4), assume $\zeta_m\subseteq F$.
If $m\equiv 2$ (mod 4), assume $\zeta_{2m}\subseteq F$.
Then $H^n(G;\mathbb{Z})_{\text{neg},F}=H^n(G;\mathbb{Z})$ for $d\in\{4,5\}$.
\end{corollary}

\begin{proof}
Since multiplication by $m$ on $H^n(G;\mathbb{Z})$ is trivial in each degree, short exact sequence
\begin{equation*}
\begin{tikzcd}
    0\arrow{r} &\mathbb{Z}\arrow{r}{m} &\mathbb{Z}\arrow{r} &\mathbb{Z}/m\mathbb{Z}\arrow{r} &0
\end{tikzcd}
\end{equation*}
induces a surjective connecting map on cohomology $H^n(G;\mathbb{Z}/m\mathbb{Z})\to H^{n+1}(G;\mathbb{Z})$ for each $d\geq 1$.
The connecting map respects negligible cohomology since the short exact sequence produces an analogous long exact sequence in Galois cohomology.
Proposition \ref{cyclic Z/mZ coefficients} completes the argument in degrees 4 and 5.

\end{proof}

Let $G$ be a cyclic group of order $m$ with $\ell|m$.
If $m\neq 2$, then the cohomology ring
\begin{equation*}
    H^*(G;\mathbb{Z}/\ell\mathbb{Z})\cong \Lambda_{\cyc{\ell}}(G^*)\otimes\mathbb{Z}/\ell\mathbb{Z}[y]
\end{equation*}
for the characters $G^*$ generated by $x$ in degree one and $\deg(y)=2$ by \cite[Section 3.2]{Ev}.
When $m=2$,
\begin{equation*}
    H^*(G;\mathbb{Z}/2\mathbb{Z})\cong\mathbb{Z}/2\mathbb{Z}[x]
\end{equation*}
for $\deg(x)=1$ by \cite[Section 3.2]{Ev}.
Furthermore,
\begin{equation*}
    H^*(G;\mathbb{Z})\cong\mathbb{Z}/n\mathbb{Z}[y]
\end{equation*}
for $\deg(y)=2$.
Since inflation maps commute with cup products, the cup product of a negligible element of $H^*(G;\mathbb{Z}/m\mathbb{Z})$ with any element of $H^*(G;\mathbb{Z}/m\mathbb{Z})$ will be negligible.

\begin{proof}[Proof of Theorem \ref{nthresh cyclic group}]
Let $\ell|m$.
The short exact sequence
\begin{equation*}
\begin{tikzcd}
    0\arrow{r} &\mathbb{Z}/\left(\frac{m}{\ell}\right)\mathbb{Z}\arrow{r} &\mathbb{Z}/m\mathbb{Z}\arrow{r} &\mathbb{Z}/\ell\mathbb{Z}\arrow{r} &0
\end{tikzcd}
\end{equation*}
induces
\begin{equation*}
\begin{tikzcd}
    0\arrow{r} &H^n\left(G;\mathbb{Z}/\left(\frac{m}{\ell}\right)\mathbb{Z}\right)\arrow{r} &H^n\left(G;\mathbb{Z}/m\mathbb{Z}\right)\arrow{r} &H^n\left(G;\mathbb{Z}/\ell\mathbb{Z}\right)\arrow{r} &0
\end{tikzcd}
\end{equation*}
in each degree $d\geq0$.
Thus, $H^n(G;\mathbb{Z}/m\mathbb{Z})\to H^n(G;\mathbb{Z}/\ell\mathbb{Z})$ is surjective for all $d\geq0$.
By Lemma \ref{nat} and Proposition \ref{cyclic Z/mZ coefficients}, we conclude that $H^3(G;\mathbb{Z}/\ell\mathbb{Z})$ and $H^4(G;\mathbb{Z}/\ell\mathbb{Z})$ are entirely negligible.
The ring structure of $H^*(G;\mathbb{Z}/\ell\mathbb{Z})$ proves that $\nthres{F}{G}{\mathbb{Z}/\ell\mathbb{Z}}\leq 3$.
By \cite[Corollary 2.2]{GM}, $2\leq\nthres{F}{G}{\mathbb{Z}/\ell\mathbb{Z}}$ as desired.

Proposition \ref{basics}(2) shows $3\leq\nthres{F}{G}{\mathbb{Z}}$. Corollary \ref{cyclic Z coefficients} and the ring structure of $H^*(G;\mathbb{Z})$ implies $\nthres{F}{G}{\mathbb{Z}}\leq4$.
Further, $H^3(G;\mathbb{Z})$ is trivial for $G$ cyclic by \cite[Section III.1]{brown1994cohomology}, completing the proof.

\end{proof}

\section{Bounds on the Negligible Thresholds Mod-$p$ of General Finite Abelian $p$-Groups}\label{7. negligible thoreslds mod-p of general finite abelian p-groups}

\subsection{The odd prime case}

Let $G\cong\bigoplus_{i=1}^m G_i$ be a finite abelian $p$-group for odd prime $p$ with $G_i\cong\mathbb{Z}/p^{k_i}\mathbb{Z}$ for $k_i\in\mathbb{N}_{>0}$.
By \cite[Sections 3.2 and 3.5]{Ev},
\begin{equation*}
    H^*(G;\mathbb{Z}/p\mathbb{Z})\cong\Lambda_{\mathbb{Z}/p\mathbb{Z}}(G^*)[y_1,\dots,y_n]
\end{equation*}
where $G^*$ is the characters of $G$ and $\deg(y_i)=2$ for $1\leq i\leq n$.
We can pick a basis $\{x_1,\dots,x_n\}$ for $G^*=H^1(G;\mathbb{Z}/p\mathbb{Z})$ as an $\mathbb{F}_p$-vector space such that $x_i$ is the inflation of a generator of $H^1(G_i;\mathbb{Z}/p\mathbb{Z})$.
Further, we can pick $y_i$ so that it is the inflation of a generator of $H^2(G_i;\mathbb{Z}/p\mathbb{Z})$.

\begin{lemma}\label{inflation1}
Let $F$ be a field such that $\text{char}(F)\neq p$ and $\zeta_{\text{exp}(G)}\subseteq F$.
The set $\{x_iy_i,y_i^2\}$ is negligible in $H^*(G;\mathbb{Z}/p\mathbb{Z})$ over $F$ for each $1\leq i\leq n$.
\end{lemma}

\begin{proof}
By Proposition \ref{cyclic Z/mZ coefficients}, $x_iy_i$ and $y_i^2$ are negligible in $H^3(G_i;\mathbb{Z}/|G_i|\mathbb{Z})$ and $H^4(G_i;\mathbb{Z}/|G_i|\mathbb{Z})$, respectively.
The quotient $\mathbb{Z}/|G_i|\mathbb{Z}\to\mathbb{Z}/p\mathbb{Z}$ induces a surjective map $H^n(G_i;\mathbb{Z}/|G_i|\mathbb{Z})\to H^n(G_i;\mathbb{Z}/p\mathbb{Z})$ for all $n\geq0$.
Thus, $x_iy_i$ and $y_i^2$ are negligible in $H^3(G_i;\mathbb{Z}/p\mathbb{Z})$ and $H^4(G_i;\mathbb{Z}/p\mathbb{Z})$, respectively.
By \cite[Proposition 2.3(a)]{GM}, the inflation map $\text{inf}:H^*(G_i;\mathbb{Z}/p\mathbb{Z})\to H^*(G;\mathbb{Z}/p\mathbb{Z})$ respects negligible classes, so $x_iy_i$ and $y_i^2$ are negligible in $H^3(G;\mathbb{Z}/p\mathbb{Z})$ and $H^4(G;\mathbb{Z}/p\mathbb{Z})$, respectively.

\end{proof}

\begin{theorem}\label{finite abelian odd p}
Let $G\cong\bigoplus_{i=1}^m\mathbb{Z}/p^{k_i}\mathbb{Z}$ be a finite abelian $p$-group for odd prime $p$ and $k_i\in\mathbb{N}^+$.
\begin{enumerate}[(1)]
    \item Let $F$ be a field such that $\zeta_{\text{exp}(G)}\subseteq F$.
    Then $\nthres{F}{G}{\mathbb{Z}/p\mathbb{Z}}\leq 2n+1$.
    \item Let $F$ be a field such that $\zeta_{p\cdot\text{exp}(G)}\subseteq F$.
    Then $\nthres{F}{G}{\mathbb{Z}/p\mathbb{Z}}\leq n+1$.
\end{enumerate}
\end{theorem}

\begin{proof}
\begin{enumerate}[(1)]
    \item Let $u\in H^d(G;\mathbb{Z}/p\mathbb{Z})$ be a non-zero monomial for $d\geq 2n+1$.
    By pigeonhole principle, $u$ is divisible by $x_iy_i$ or $y_i^2$ for some $1\leq i\leq n$.
    Thus, $u$ is negligible over $F$ by Lemma \ref{inflation1}.
    We conclude that every element of $H^d(G;\mathbb{Z}/p\mathbb{Z})$ is negligible for $d\geq 2n+1$.
    \item Let $u\in H^d(G;\mathbb{Z}/p\mathbb{Z})$ for $d\geq n+1$.
    By pigeonhole principle, $u$ is divisible by $y_i$ for some $1\leq i\leq n$.
    In the cohomology of the $i$th summand of $G$, \cite[Corollary 3.4]{GM} proves that the class that corresponds to $y_i$ is negligible.
    By \cite[Proposition 2.3(a)]{GM}, the inflation map respects negligible classes, which implies that $y_i$ is negligible in $H^2(G;\mathbb{Z}/p\mathbb{Z})$.
    We conclude that every element of $H^d(G;\mathbb{Z}/p\mathbb{Z})$ is negligible for $d\geq n+1$.
\end{enumerate}

\end{proof}

\begin{remark}
Let $p$ be an odd prime and $G$ a cyclic $p$-group.
Let $F$ be a field such that $\zeta_{|G|}\subseteq F$ but $\zeta_{p|G|}\not\subseteq F$.
The bound in Theorem \ref{finite abelian odd p}(2) is sharp for $\nthres{F}{G}{\mathbb{Z}/p\mathbb{Z}}$ by Theorem \ref{nthresh cyclic group} and \cite[Corollary 3.4]{GM}.
Let $E$ be an elementary abelian $p$-group for $p$ odd.
The bound in Theorem \ref{finite abelian odd p}(3) is sharp for $\nthres{F}{E}{\mathbb{Z}/p\mathbb{Z}}$ by Theorem \ref{elementary abelian threshold} since we implicitly assume $\zeta_{p^2}\subseteq F$.
The bound in Theorem \ref{finite abelian odd p}(1) is probably not sharp.
\end{remark}

\subsection{The case $p=2$}

Let $G$ be a finite abelian 2-group.
Assume that $G$ has $n$ cyclic direct summands in elementary divisor form.
Denote by $G_i$ the $i$th cyclic direct summand of $G$ so $G\cong\bigoplus_{i=1}^nG_i$.
Denote by $m$ the first index for which $|G_i|>2$.
By \cite[Sections 3.2 and 3.5]{Ev}, the cohomology ring of $A$ is 
\begin{equation*}
    H^*(G;\mathbb{Z}/2\mathbb{Z})\cong\mathbb{Z}/2\mathbb{Z}[x_1,\dots,x_n,y_m,\dots,y_n]/\langle x_m^2,\dots,x_n^2\rangle
\end{equation*}
where $\deg(x_i)=1$ for $1\leq i\leq n$ and $\deg(y_j)=2$ for $m\leq j\leq n$.
We can pick a basis $\{x_1,\dots,x_n\}$ for $H^1(G;\mathbb{Z}/2\mathbb{Z})$ as a $\mathbb{Z}/2\mathbb{Z}$-vector space such that $x_i$ is the inflation of a generator of $H^1(G_i;\mathbb{Z}/2\mathbb{Z})$.
Further, we can pick $y_i$ so that it is the inflation of a generator of $H^2(G_i;\mathbb{Z}/2\mathbb{Z})$.

\begin{lemma}\label{inflation2}
Let $\zeta_{4}\subseteq F$.
The elements $x_i^2$ for each $1\leq i\leq m-1$ and $y_j^2$ for each $m\leq j\leq n$ are negligible in $H^*(G;\mathbb{Z}/2\mathbb{Z})$ over $F$.
\end{lemma}

\begin{proof}
Let $K$ be a field extension of $F$ with group homomorphism $j:\Gamma_K\to G_i$.
Let $x_i'$ be a generator of $H^1(G_i;\mathbb{Z}/2\mathbb{Z})$ for $1\leq i\leq m-1$ and $y_j'$ be a generator of $H^2(G_j;\mathbb{Z}/2\mathbb{Z})$ for $m\leq j\leq n$.
Then $j^*((x_i')^2)=j^*(x_i')\cup(-1)=0$ and $j^*((y_j')^2)=j^*(y_j')\cup(-1)^2=0$ by the Norm Residue Isomorphism Theorem \cite{HW}.
We conclude $(x_i')^2$ and $(y_j')^2$ are negligible over $F$.
By \cite[Proposition 2.3(a)]{GM}, the inflation map $\text{inf}:H^*(G_i;\mathbb{Z}/2\mathbb{Z})\to H^*(G;\mathbb{Z}/2\mathbb{Z})$ respects negligible classes.

\end{proof}

\begin{lemma}\label{inflation3}
Let $\zeta_{2\text{exp}(G)}\subseteq F$.
The elements $x_i^2$ for each $1\leq i\leq m-1$ and $y_j$ for each $m\leq j\leq n$ are negligible in $H^*(G;\mathbb{Z}/2\mathbb{Z})$ over $F$.
\end{lemma}

\begin{proof}
A computation using \cite[Corollary 3.4]{GM} proves that generator $(x_i')^2$ is negligible in $H^2(G_i,\mathbb{Z}/p\mathbb{Z})$ and generator $y_j'$ is negligible in $H^2(G_j;\mathbb{Z}/p\mathbb{Z})$.
By \cite[Proposition 2.3(a)]{GM}, the inflation map $\text{inf}:H^*(G_i;\mathbb{Z}/p\mathbb{Z})\to H^*(G;\mathbb{Z}/p\mathbb{Z})$ respects negligible classes.

\end{proof}

\begin{theorem}\label{finite abelian p=2}
Let $G\cong\bigoplus_{i=1}^m\mathbb{Z}/2^{k_i}\mathbb{Z}$ be a finite abelian 2-group for $k_i\in\mathbb{N}^+$.
\begin{enumerate}[(1)]
    \item Let $F$ be a field such that $\zeta_{4}\subseteq F$.
    Then $\nthres{F}{G}{\mathbb{Z}/2\mathbb{Z}}\leq 3n+1$.
    \item Let $F$ be a field such that $\zeta_{2\text{exp}(G)}\subseteq F$.
    Then $\nthres{F}{G}{\mathbb{Z}/2\mathbb{Z}}\leq n+1$.
\end{enumerate}
\end{theorem}

\begin{proof}
\begin{enumerate}[(1)]
    \item Let $u\in H^d(G;\mathbb{Z}/2\mathbb{Z})$ be a non-zero monomial for $d\geq 3n+1$.
    By pigeonhole principle, $u$ is divisible by $x_i^2$ for some $1\leq i\leq m-1$ or $y_j^2$ for some $m\leq j\leq n$.
    Thus, $u$ is negligible over $F$ by Lemma \ref{inflation2}.
    We conclude that every element of $H^d(G;\mathbb{Z}/2\mathbb{Z})$ is negligible for $d\geq 3n+1$.
    \item Let $u\in H^d(G;\mathbb{Z}/2\mathbb{Z})$ be a non-zero monomial for $d\geq n+1$.
    By pigeonhole principle, $u$ is divisible by $x_i^2$ or $y_j$ for some $1\leq i\leq m-1$ and $m\leq j\leq n$.
    Thus, $u$ is negligible over $F$ by Lemma \ref{inflation3}.
    We conclude that every element of $H^d(G;\mathbb{Z}/2\mathbb{Z})$ is negligible for $d\geq n+1$.
\end{enumerate}
\end{proof}

\begin{remark}
Let $G$ a cyclic 2-group.
Let $F$ be a field such that $\zeta_{2|G|}\subseteq F$.
The bound in Theorem \ref{finite abelian p=2}(2) is sharp for $\nthres{F}{G}{\mathbb{Z}/p\mathbb{Z}}$ by \cite[Corollary 3.4]{GM} and the ring structure of $H^*(G;\mathbb{Z}/2\mathbb{Z})$.
Let $E$ be an elementary abelian 2-group.
The bound in Theorem \ref{finite abelian p=2}(2) is also sharp for $\nthres{F}{G}{\mathbb{Z}/2\mathbb{Z}}$ by Theorem \ref{elementary abelian 2-group thresh} since we implicitly assume $\zeta_4\subseteq F$.
The bound in Theorem \ref{finite abelian p=2}(1) is probably not sharp.
\end{remark}

















\newpage
\bibliographystyle{plain}
\bibliography{negthres}

@book{lang_algebra_2005,
	address = {New York},
	edition = {Rev. 3. ed. (corr. print. 2005), [Nachdr.]},
	series = {Graduate texts in mathematics},
	title = {Algebra},
	isbn = {978-0-387-95385-4},
	language = {eng},
	number = {211},
	publisher = {Springer},
	author = {Lang, Serge},
	year = {2005},
}

@article{gherman_krull,
	title = {Krull dimension of the negligible quotient in mod p cohomology of a finite group},
	volume = {228},
	issn = {00224049},
	url = {https://linkinghub.elsevier.com/retrieve/pii/S0022404923001718},
	doi = {10.1016/j.jpaa.2023.107489},
	pages = {107489},
	number = {3},
	journaltitle = {Journal of Pure and Applied Algebra},
	shortjournal = {Journal of Pure and Applied Algebra},
	author = {Gherman, M. and Merkurjev, A.},
	urldate = {2025-12-29},
	date = {2024-03},
	langid = {english},
}

@book{weibelhomalg,
	location = {Cambridge},
	edition = {Reprint. 1997, transf. to digital print},
	title = {An introduction to homological algebra},
	isbn = {978-0-521-43500-0 978-0-521-55987-4},
	series = {Cambridge studies in advanced mathematics},
	pagetotal = {450},
	number = {38},
	publisher = {Cambridge Univ. Press},
	author = {Weibel, Charles A.},
	date = {2003},
}

@book{cigc,
	location = {Providence, Rhode Island},
	title = {Cohomological invariants in Galois cohomology},
	isbn = {978-0-8218-3287-5 978-1-4704-2174-8},
	series = {University lecture series},
	pagetotal = {1},
	number = {Volume 28},
	publisher = {American Mathematical Society},
	author = {Garibaldi, S. and Merkurjev, A. and Serre, J.-P.},
	date = {2003},
}

@mastersthesis{hoegenhaven2013,
  author       = {Amalie H{\o}genhaven},
  title        = {Quillen Stratification in Group Cohomology},
  school       = {K{\o}benhavns Universitet},
  year         = {2013},
  type         = {Master's thesis},
  url          = {https://web.math.ku.dk/~jg/students/hoegenhaven.msthesis.2013.pdf}
}

@incollection{BassTate1973Milnor,
  author    = {Bass, Hyman and Tate, John},
  title     = {The Milnor ring of a global field},
  booktitle = {Algebraic K-Theory, II: {\"{A}}rithmetic and Geometry},
  series    = {Lecture Notes in Mathematics},
  volume    = {342},
  pages     = {349--446},
  publisher = {Springer-Verlag},
  year      = {1973},
  address   = {Berlin}
}

@book{totaro_group_2014,
	location = {Cambridge},
	title = {Group cohomology and algebraic cycles},
	isbn = {978-1-139-05948-0},
	series = {Cambridge tracts in mathematics},
	pagetotal = {1},
	number = {204},
	publisher = {Cambridge University Press},
	author = {Totaro, Burt},
	date = {2014},
	doi = {10.1017/CBO9781139059480},
}

@book{HW,
  title={The Norm Residue Theorem in Motivic Cohomology:(AMS-200)},
  author={Haesemeyer, C. and Weibel, C.},
  year={2019},
  publisher={Princeton University Press}
}

@book{Ev,
title={The Cohomology of Groups},
author={Evens, L.},
year={1991},
publisher={Clarendon Press, Oxford}
}

@article {GM,
    AUTHOR = {Gherman, M. and Merkurjev, A.},
     TITLE = {Negligible degree two cohomology of finite groups},
   JOURNAL = {J. Algebra},
  FJOURNAL = {Journal of Algebra},
    VOLUME = {611},
      YEAR = {2022},
     PAGES = {82--93},
      ISSN = {0021-8693},
   MRCLASS = {20J06 (12F10)},
  MRNUMBER = {4472425},
MRREVIEWER = {Yuval Ginosar},
       DOI = {10.1016/j.jalgebra.2022.07.039},
       URL = {https://doi.org/10.1016/j.jalgebra.2022.07.039},
}

@book {Se,
    AUTHOR = {Serre, J.-P.},
     TITLE = {Oeuvres/{C}ollected papers. {IV}. 1985--1998},
    SERIES = {Springer Collected Works in Mathematics},
      NOTE = {Reprint of the 2000 edition [MR1730973]},
 PUBLISHER = {Springer, Heidelberg},
      YEAR = {2013},
     PAGES = {viii+694},
      ISBN = {978-3-642-39839-1},
   MRCLASS = {01A75 (01A70)},
  MRNUMBER = {3185222},
}

@article{MS24,
    title={Galois representations modulo $p$ that do not lift modulo $p^2$},
    author={A. Merkurjev and F. Scavia},
    url={https://arxiv.org/abs/2410.12560}
}

@misc{KL20,
      title={Liftable groups, negligible cohomology and Heisenberg representations}, 
      author={Chandrashekhar B. Khare and Michael Larsen},
      year={2020},
      eprint={2009.01301},
      archivePrefix={arXiv},
      primaryClass={math.NT},
      url={https://arxiv.org/abs/2009.01301}, 
}

@book{elman2008algebraic,
  title={The Algebraic and Geometric Theory of Quadratic Forms},
  author={Elman, R.S. and Karpenko, N. and Merkurjev, A.},
  isbn={9780821873229},
  series={American Mathematical Society colloquium publications},
  url={https://books.google.com/books?id=NsYPMC3iBowC},
  year={2008},
  publisher={American Mathematical Soc.}
}

@book{brown1994cohomology,
  title={Cohomology of Groups},
  author={Brown, K.S.},
  isbn={9780387906881},
  lccn={82000733},
  series={Graduate Texts in Mathematics},
  url={https://books.google.com/books?id=PMqb2DppvCsC},
  year={1994},
  publisher={Springer New York}
}

\end{document}